\documentclass[12pt]{amsart}

\usepackage{latexsym}
\usepackage{psfrag}
\usepackage{amsmath}
\usepackage{amssymb}
\usepackage{epsfig}
\def\LOCAL{.}
\usepackage{amsfonts}

\def\unit#1{\refstepcounter{section}\part{#1}}

\newcommand{\excise}[1]{}

\newtheorem{thm}{Theorem}[section]
\newtheorem{lemma}[thm]{Lemma}
\newtheorem{claim}[thm]{Claim}
\newtheorem{cor}[thm]{Corollary}
\newtheorem{prop}[thm]{Proposition}
\newtheorem{conj}[thm]{Conjecture}
\newtheorem{conv}[thm]{Convention}

\newtheorem{Example}[thm]{Example}
\newtheorem{Remark}[thm]{Remark}
\newtheorem{Defn}[thm]{Definition}
\newtheorem{Alg}[thm]{Algorithm}
\newtheorem{Rtne}[thm]{Routine}
\newtheorem{Warn}[thm]{Caution}

\newenvironment{example}{\begin{Example}\rm}
    {\end{Example}}
\newenvironment{remark}{\begin{Remark}\rm}
    {\end{Remark}}
\newenvironment{defn}{\begin{Defn}\rm}{\end{Defn}}
\newenvironment{alg}{\begin{Alg}\rm}{\end{Alg}}
\newenvironment{rtne}{\begin{Rtne}\rm}{\end{Rtne}}

\newenvironment{proofof}[1]{\begin{trivlist}\item {\bf
    Proof of {#1}.\,}}{\mbox{}\hfill$\square$\end{trivlist}}

\newcounter{separated}

\def\zz{\mathbb Z}

\def\rr{\mathbb R}
\def\qqq{\mathbb Q}

\def\pr{\prime}
\def\Ga{\Gamma}

\def\la{\lambda}
\def\ga{\gamma}

\def\al{\alpha}

\def\om{\omega}

\def\ve{\varepsilon}

\def\vp{\varphi}

\def\CR{\mathcal R}

\def\cf{\mathcal F}

\def\cl{\mathcal L}
\def\co{\mathcal O}

\def\ssu{\subset}

\def\<{\langle}
\def\>{\rangle}

\def\Ups{\Upsilon}

\def\0{{\mathbf 0}}

\def\st{{\rm st}}
\def\th{{\rm th}}
\def\ol#1{{\overline {#1}}}

\def\EE{{\mathcal E}}
\def\II{{\mathcal I}}
\def\Kv{K_{\hspace{-.2ex}v}}
\def\KK{\hspace{.35ex}\ol{\hspace{-.35ex}K\hspace{-.05ex}}\hspace{.05ex}}
\def\KKv{\KK_{\hspace{-.2ex}v}}
\def\LL{{\mathcal L}}
\def\NN{{\mathbb N}}
\def\OO{{\mathcal O}}
\def\TT{{\mathcal T}}
\def\VV{{\mathcal V}}

\def\TF{T_{\hspace{-.15ex}F}}
\def\TFp{T_{\hspace{-.15ex}F'}}
\def\TFn#1{T_{\hspace{-.15ex}F_{#1}}}
\def\TR{T_{\hspace{-.2ex}R}}

\def\Tv{T_{\hspace{-.12ex}v}}
\def\Tw{T_{\hspace{-.12ex}w}}
\def\Tx{T_{\hspace{-.12ex}x}}

\def\Txi{T_{\hspace{-.12ex}\xi}}
\def\Txx#1{T_{\hspace{-.12ex}x}^{\hspace{.12ex}#1\hspace{.12ex}}}
\def\UF{\Upsilon_{\!F}}
\def\UFp{\Upsilon_{\!F'}}
\def\UG{\Upsilon_{\!G}}

\def\UUv{\ol U_{\hspace{-.25ex}v}}
\def\Uv{U_{\hspace{-.12ex}v}}

\def\hce{\makebox[0ex]{\raisebox{.4ex}{\hspace{1.7ex}$\widehat{}$}}\mathcal E}
\def\src{{\rm src}}

\def\minus{\smallsetminus}

\def\upper{\hspace{.3ex}\overline{\mathit{src}}_v}
\def\uppers{\hspace{.3ex}\overline{\mathit{src}}}

\def\implies{\Rightarrow}
\def\nothing{\varnothing}

\def\EEE{\hce}

\def\bot{\text{bot}}
\def\top{\text{top}}
\def\back{\text{back}}
\def\fro{\text{front}}
\def\rig{\text{right}}
\def\lef{\text{left}}
\newenvironment{alglist}%
    {\begin{list}
        {}
        {\leftmargin=5.4em\labelwidth=5.1em\labelsep=.6em
         \topsep=-1ex\itemsep=.1ex}\sf}  
    {\vspace{1ex}\end{list}}
\def\routine#1{\item[{\sc{#1}{\ }}]}
\def\procedure#1{{\sc{#1}}}
\newenvironment{routinelist}[1]%
    {\routine{#1}\begin{list}
        {}
        {\leftmargin=2.9em\labelwidth=2.4em\labelsep=.5em
         \topsep=-1ex\itemsep=.1ex}}
    {\end{list}}
    {\begin{list}
        {}
        {\leftmargin=3.0em\labelwidth=4.8em\labelsep=.5em
         \itemsep=.1ex\topsep=0ex}}
    {\end{list}}

\newenvironment{numbered}%
        {\begin{list}
                {\noindent\makebox[0mm][r]{\arabic{enumi}.}}
                {\leftmargin=5.5ex \usecounter{enumi}}
        }
        {\end{list}}

        {\begin{list}
                {\noindent\makebox[0mm][r]{(\roman{enumi})}}
                {\leftmargin=5.5ex \usecounter{enumi}}
        }
        {\end{list}}

\begin{document}

\title[Metric combinatorics of convex polyhedra]%
    {Metric combinatorics of convex polyhedra:\\ Cut loci and
    nonoverlapping unfoldings}
\author[Ezra Miller]{Ezra Miller$^*$}
\author[Igor~Pak]{Igor~Pak$^\dagger$}
\date{4 December 2003}


\thanks{${\hspace{-1ex}}^*$School of Mathematics, University of
Minnesota, Minneapolis, MN;~\texttt{ezra@math.umn.edu}}

\thanks{${\hspace{-.95ex}}^\dagger$Department of Mathematics,
MIT, Cambridge, MA;~\texttt{pak@math.mit.edu}}

\begin{abstract}
\noindent
Let~$S$ be the boundary of a convex polytope of
dimension~\mbox{$d+1$}, or more generally let~$S$ be a {\em convex
polyhedral pseudomanifold}.  We prove that~$S$ has a polyhedral
nonoverlapping unfolding into~$\rr^d$, so the metric space $S$ is
obtained from a closed (usually nonconvex) polyhedral ball
in~$\rr^d$ by identifying pairs of boundary faces isometrically.
Our existence proof exploits geodesic flow away from a source point
$v \in S$, which is the exponential map to~$S$ from the tangent
space at~$v$.  We characterize the \emph{cut locus} (the closure of the
set of points in~$S$ with more than one shortest path to~$v$) as a
polyhedral complex in terms of Voronoi diagrams on facets.
Analyzing infinitesimal expansion of the wavefront consisting of
points at constant distance from~$v$ on~$S$ produces an algorithmic
method for constructing Voronoi diagrams in each facet, and hence
the unfolding of~$S$.  The algorithm, for which we provide
pseudocode, solves the discrete geodesic problem.  Its main
construction generalizes the source unfolding for boundaries of
$3$-polytopes into~$\rr^2$.  We present conjectures concerning the
number of shortest paths on the boundaries of convex polyhedra, and
concerning continuous unfolding of convex polyhedra.  We also
comment on the intrinsic non-polynomial complexity of nonconvex
manifolds.
\end{abstract}

\maketitle


\part*{Contents}

\begin{list}{\arabic{enumi}.}
           {\leftmargin=5ex \rightmargin=2ex \usecounter{enumi}}
\item[{\makebox[1.7ex][l]{Introduction}}]
\hfill \pageref{intropage}
\item[{\makebox[1.7ex][l]{Overview}}]
\hfill \pageref{overviewpage}
\item[{\makebox[1.7ex][l]{Methods}}]
\hfill \pageref{methodspage}
\item
Geodesics in polyhedral boundaries \hfill \pageref{geodesics}
\item
Cut loci \hfill \pageref{cutloci}
\item
Polyhedral nonoverlapping unfolding \hfill \pageref{unfold}
\item
The source poset \hfill \pageref{wavefront}
\item
Constructing source images \hfill \pageref{construct}
\item
Algorithm for source unfolding \hfill \pageref{algorithm}
\item
Convex polyhedral pseudomanifolds \hfill \pageref{pseudomanifolds}
\item
Limitations, generalizations, and history \hfill \pageref{limits}
\item
Open problems and complexity issues \hfill \pageref{open}
\item[{\makebox[1.7ex][l]{Acknowledgments}}]
\hfill \pageref{ackpage}
\item[{\makebox[1.7ex][l]{References}}]
\hfill \pageref{refpage}
\end{list}

\part*{Introduction}\label{intropage}

\noindent
The past several decades have seen intense development in the
combinatorics and geometry of convex polytopes \cite{Z}.  Besides
their intrinsic interest, the advances have been driven by
applications to areas ranging as widely as combinatorial
optimization, commutative algebra, symplectic geometry, theoretical
physics, representation theory, statistics, and enumerative combinatorics.
As a result, there is currently available a wealth of
insight into (for example) algebraic invariants of the face posets
of polytopes;
arithmetic information connected to sets of lattice points inside
polytopes;
and geometric constructions associated with linear functionals,
such as Morse-like decompositions and methods for locating extrema.

On the topological side, there are metric theories for polyhedral
spaces, primarily motivated by differential geometry.
In addition, there is a vast literature on general convexity.
Nonetheless, there seems to be lacking a study of the interaction
between the combinatorics of the boundaries of convex polytopes and
their metric geometry in arbitrary dimension.  This remains the
case despite relations to a number of classical algorithmic
problems in discrete and computational geometry.

The realization here is that convexity and polyhedrality together
impose rich combinatorial structures on the collection of shortest
paths in a metrized sphere.  We initiate a systematic investigation
of this {\em metric combinatorics}\/ of convex polyhedra by proving
the existence of polyhedral nonoverlapping unfoldings and analyzing
the structure of the cut locus.  The algorithmic aspect, which we
include together with its complexity analysis, was for us a
motivating feature of these results.  That being said, we also show
that our general methods are robust enough so that---with a few
minor modifications---they extend to the abstract spaces we call
`convex polyhedral pseudomanifolds', whose sectional curvatures
along low-dimensional faces are all positive.  To conclude, we
propose some directions for future research, including a series of
precise conjectures on the number of combinatorial types of
shortest paths, and on the geometry of unfolding boundaries of
polyhedra.

\vspace{-1ex}
\part*{Overview}\label{overviewpage}

\noindent
Broadly speaking, the metric geometry of boundaries of
$3$-dimensional polytopes is quite well understood, due in large
part to work of A.\thinspace{}D.\thinspace{}Aleksandrov
\cite{Al48,Al50} and his school.  For higher dimensions, however,
less theory appears in the literature, partly because Aleksandrov's
strongest methods do not extend to higher dimension.  Although
there do exist general frameworks for dealing with metric geometry
in spaces general enough to include boundaries of convex polyhedra,
such as \cite{BGP}, the special nature of polyhedral spaces usually
plays no role.

The existing theory that does appear for polyhedral spaces is
motivated from the perspective of Riemannian geometry, via metric
geometry on simplicial complexes, and seems mainly due to
D.\thinspace{}A.\thinspace{}Stone; see \cite{Sto76}, for example.
In contrast, our original motivation comes from two classical
problems in discrete and computational geometry: the `discrete
geodesic problem'~\cite{Mit} of finding shortest paths between
points on polyhedral surfaces, and the problem of constructing
nonoverlapping unfoldings of convex polytopes~\cite{O}.  Both
problems are well understood for the $2$-dimensional boundaries of
$3$-polytopes, but have not been attempted in higher dimensions.
We resolve them here in arbitrary dimension by a unified
construction generalizing the `source unfolding' of
\mbox{$3$-dimensional} convex polyhedra~\cite{VP71,SS}.

Previous methods for source unfoldings have been specific to low
dimension, relying for example on the fact that arcs of circles in the
plane intersect polygons in finite sets of points.  We instead use
techniques based on differential geometry to obtain general results
concerning cut loci on boundaries of polytopes in arbitrary dimension,
namely Theorem~\ref{vor} and Corollary~\ref{c:exp}, thereby producing
polyhedral foldouts in~Theorem~\ref{exp}.  In more precise terms, our
two main goals in this paper~are~to:
\begin{numbered}
\item
describe how the set of points on the boundary~$S$ of a convex
polyhedron at given radius from a fixed {\em source point}\/ changes
as the radius increases continuously;

\item
use this description of `wavefront expansion' to construct a
polyhedral nonoverlapping unfolding of the $d$-dimensional polyhedral
complex~$S$ into~$\rr^d$.
\end{numbered}
By `describe' and `construct' we mean to achieve these goals not
just abstractly and combinatorially, but effectively, in a manner
amenable to algorithmic computation.  References such as
\cite{AAOS,AO,CH,Ka,MMP,Mount85,SS}, which have their roots and
applications
%
%
in computational geometry, carry this out in the \mbox{$d=2$} case
of boundaries of $3$-polytopes (and for the first goal, on any
polyhedral surface of dimension $d=2$).  Here, in arbitrary
dimension~$d$, our Theorem~\ref{src} says precisely how past
wavefront evolution determines the location in time and space of
its next qualitative change.  The combinatorial nature of
Theorem~\ref{src} leads immediately to Algorithm~\ref{code} for
effectively unfolding boundaries of polyhedra.

The results and proofs in Sections~\ref{geodesics}--\ref{algorithm}
for boundaries of convex polyhedra almost all hold verbatim in the
more abstract setting of what we call $d$-dimensional {\em convex
polyhedral pseudomanifolds}.  The study of such spaces is suggested
both by Stone's point of view in \cite{Sto76} and by the more
general methods in \cite{BGP}.  Our Corollary~\ref{star-shaped}
says that all convex polyhedral pseudomanifolds can be represented
as quotients of Euclidean (usually nonconvex) polyhedral balls by
identifying pairs of boundary components isometrically.  The reader
interested solely in this level of generality is urged to begin
with Section~\ref{pseudomanifolds}, which gives a guide to
Sections~\ref{geodesics}--\ref{algorithm} from that perspective,
and provides the slight requisite modifications where necessary.
Hence the reader can avoid checking the proofs in the earlier
sections twice.

The results in Section~\ref{pseudomanifolds} on convex polyhedral
pseudomanifolds are in many senses sharp, in that considering more
general spaces would falsify certain conclusions.  We substantiate
this claim in Section~\ref{limits}, where we also discuss
extensions of our methods that are nonetheless possible.  For
example, we present an algorithm to construct geodesic Voronoi
diagrams on boundaries of convex polyhedra in
Section~\ref{geod-voronoi}.

The methods of this paper suggest a number of fundamental open
questions about the metric combinatorics of convex polyhedra in
arbitrary dimension, and we present these in Section~\ref{open}.
Most of them concern the notion of {\em vistal tree}\/ in
Definition~\ref{tree}, which encodes all of the combinatorial types
of shortest paths (or equivalently, all bifurcations of the
wavefront) emanating from a source point.  The first two questions,
Conjectures~\ref{polynomial} and~\ref{polynomial'}, concern the
complexity of our unfolding algorithm and the behavior of geodesics
in boundaries of polyhedra.  Along these lines, we remark also on
the complexity of nonconvex polyhedral manifolds, in
Proposition~\ref{nonconvex}.  Our third question is about the
canonical subdivision of the boundary of any convex polyhedron
determined by the sets of source points having isomorphic vistal
trees (Definition~\ref{equivistal} and Conjecture~\ref{vistal}); it
asks whether this {\em vistal subdivision}\/ is polyhedral, and how
many faces it has.  Our final question asks how to realize
unfoldings of polyhedral boundaries by embedded homotopies
(Conjecture~\ref{bloom}).

\vspace{-1ex}
\part*{Methods}\label{methodspage}

This section contains an extended overview of the paper, including
background and somewhat informal descriptions of the geometric
concepts involved.

\subsection*{Unfolding polyhedra}

While unfolding convex polytopes is easy~\cite{Al48}, constructing
a \emph{nonoverlapping}\/ unfolding is in fact a difficult task
with a long history going back to D\"urer in 1528~\cite{O}.  When
cuts are restricted to ridges (faces of dimension~\mbox{$d-1$} in a
polyhedron of dimension~\mbox{$d+1$}), the existence of such
unfoldings is open even for polytopes in~$\rr^3$~\cite{O,Z}.  It is
known that nonconvex polyhedral surfaces need not admit such
nonoverlapping unfoldings~\cite{BDEK}.

In this paper we consider unfoldings of a more general
nature---when cuts are allowed to slice the interiors of facets.
Such unfoldings are known, but only for three dimensional
polytopes~\cite{SS,AO,CH,AAOS}.  In fact, two different (although
strongly related) unfoldings appear in these and other references
in the literature: the {\em Aleksandrov unfolding}\/ (also known as
the {\em star unfolding}) \cite{Al48,AO}, and the {\em source
unfolding}\/~\cite{VP71,SS}.  Unfortunately, the construction of
Aleksandrov unfoldings fails in principle in higher dimension
(Section~\ref{aleksandrov}).  As we mentioned earlier, we
generalize the source unfolding construction to prove that the
boundary~$S$ of any convex polyhedron of dimension~\mbox{$d+1$},
and more abstractly any convex polyhedral pseudomanifold~$S$, has a
nonoverlapping polyhedral unfolding~$\ol U$ in~$\rr^d$.  The second
of the following two foldouts of the cube is a $d=2$ example of a
source unfolding.
\begin{figure}[hbt]
\begin{center}
\psfrag{v}{\footnotesize $\hspace{.15ex}v$}
\epsfig{file=\LOCAL/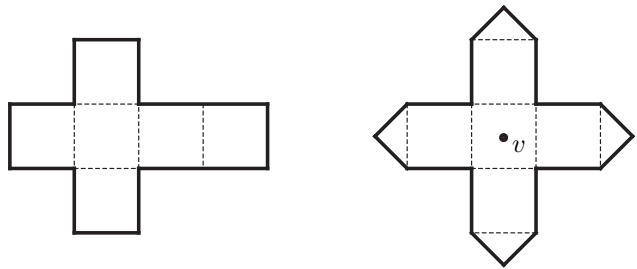,width=8.6cm}
\end{center}
\caption{An edge-unfolding and a source unfolding of a cube into
$\rr^2$}
\label{f:two-fold}
\end{figure}
For clarity, we present the discussion below in the context of
boundaries of polyhedra.

\subsection*{Cut loci}

The idea of the source unfolding in arbitrary dimension~$d$ is
unchanged from the case $d=2$ of convex polyhedral surfaces.  Pick
a {\em source point}\/~$v$ interior to some facet ($d$-dimensional
face) of~$S$, so the tangent space~$\Tv$ is well-defined.  Then,
treating~$S$ like a Riemannian manifold, define the {\em
exponential map} from $\Tv$ to~$S$ by flowing along geodesics
emanating from~$v$.  Our main unfolding result, Theorem~\ref{exp},
says that exponentiation takes a certain open polyhedral ball $\Uv
\subset \Tv$ isometrically to a dense open subset of~$S$ consisting
of points possessing a unique shortest path (length-minimizing
geodesic) to~$v$.  The image of the closure~$\UUv$ of the open
ball~$\Uv$ is all of~$S$.  The boundary $\UUv \minus \Uv$ maps onto
the {\em cut locus}~$\KKv$, which by definition is the closure of
the set of points in~$S$ with more than one shortest path to
the~source~point~$v$.  These properties characterize~$\Uv$.

In Riemannian geometry, when the manifold and the metric are both
smooth, describing the cut locus for a source point is already an
important and interesting problem (see \cite{Kob89} for an
excellent introduction and numerous references), although of course
the exponential map can only be an isometry, even locally, if the
metric is flat.  Extending the notion of cut locus from Riemannian
geometry to the polyhedral context is just as easy as extending it
to arbitrary metric spaces.  But showing that the open ball~$\Uv$
is a polyhedral foldout requires strong conditions on the
complement of the cut locus, such as metric flatness and
polyhedrality.  We prove these results in Sections~\ref{geodesics}
and~\ref{cutloci} using methods based on the foundations of
polyhedral geometry, and on Voronoi diagrams, culminating in
Theorem~\ref{vor} and Corollary~\ref{c:exp}.  These conclusions
depend crucially on convexity and do not hold in the nonconvex
case.

\subsection*{Geometry of wavefront expansion}

Our existence proof for polyhedral nonoverlapping source foldouts,
even given their Voronoi characterization in Theorem~\ref{vor},
does not by itself provide a satisfactory combinatorial picture of
the dynamics of wavefront expansion on polyhedra.  For this, we
must gain control over how the exponential map behaves as it
interacts with {\em warped}\/ points in~$S$, namely those of
infinite curvature, or equivalently points on faces of
dimension~\mbox{$d-2$} or less.

Imagine the picture kinetically: the source point~$v$ emits a signal,
whose wavefront proceeds as a $(d-1)$-sphere of increasing radius---at
least until the sphere hits the boundary of the facet containing~$v$.
At that stage, the wavefront folds over a {\em ridge}, or face of
dimension~\mbox{$d-1$}.  Metrically, nothing has happened: points
interior to ridges look to the wavefront just as flat as points
interior to facets.  But later, as the wavefront encounters faces of
lower dimension, it is forced to bifurcate around warped points and
interfere with itself, as signals emitted originally in different
directions from~$v$ begin to curl around the infinite curvature and
converge~toward~the~cut~locus.

The question becomes: What discrete structure governs evolution of
the wavefront on polyhedra?  The most obvious first step is to
define a finite collection of `events', representing the points in
time and space where the wavefront changes in some nontrivial way.
If this is done properly, then it remains only to order the events
according to the times at which they occur.  However, in reality,
the definition of an event is rather simple, while the geometry
dictating time order of events is more complex.

Starting from scratch, one might be tempted (and we were) to mark an
event every time the wavefront encounters a new warped face.  Indeed,
this works in dimension $d=2$ \cite{MMP}: since the wavefront is a
curve, its intersection with the set of edges is a finite set, and it
is easy to detect when one of these points hits a vertex of~$S$ before
another.  But because the geometry is substantially more complicated
in higher dimensions, in the end we found it more natural to say an
{\em event}\/ has occurred every time the wavefront encounters a new
facet through the relative interior of a {\em ridge} (see
Definitions~\ref{d:src} and~\ref{poset}).  This may seem
counterintuitive, since the wavefront only interacts with and curls
around faces of smaller dimension.  However, wavefront collisions with
warped points lead to intersections with ridge interiors
infinitesimally afterward.  In other words, the closest point ({\em
event point}\/) on a facet to the source point~$v$ need not lie
interior to a ridge, but can just as easily be~warped.

Again think kinetically: once the wavefront has hit a new face (of
small dimension, say), it begins to creep up each of the ridges
containing that face.  Although in a macroscopic sense the
wavefront hits all of these ridges simultaneously, it creeps up
their interiors at varying rates.  Therefore the wavefront hits
some of these ridges before others in an infinitesimal sense.  The
moral is that if one wants to detect curling of the wavefront
around warped faces, it is simpler to detect the wake of this
interaction infinitesimally on the interiors of neighboring ridges.
Sufficiently refined tangent data along ridges then discretizes the
finite set of events, thereby producing the desired `metric
combinatorics' of wavefront expansion.

\subsection*{Source poset}

Making the above moral precise occupies Section~\ref{wavefront}.
To~single out a ridge whose interior is engulfed by the wavefront
at a maximal rate (thereby making it {\em closer}\/ to the source
point) essentially is to find a ridge whose angle with the
corresponding signal ray emitted from the source is minimal.  When
$d=2$, this~means that we do not simply observe two signals hitting
vertices simultaneously, but we notice also the angles at which
they hit the edges containing those vertices.  The~edge forming the
smallest angle with its signal ray is the earlier event,
infinitesimally~beating out other potential events.  (That each
angle must be measured inside some ambient facet is just one of the
subtleties that we gloss over for now.)

To distinguish events in time macroscopically, only radii (distances
from the source) are required.  When $d = 2$, as we have just seen, a
first derivative is enough to distinguish events infinitesimally.
Generally, in dimension~\mbox{$d \geq 2$}, one needs derivatives of
order less than~$d$, or more precisely, a directional derivative
successively along each of $d-1$ orthogonal directions inside a ridge.
In Section~\ref{wavefront}, these derivatives are encoded not in
single angles, but in {\em angle sequences}
(Definition~\ref{minimal}), which provide quantitative information
about the goniometry of intersections between signal rays and the
faces of varying dimension they encounter.  More qualitative---and
much more refined---data is carried by {\em minimal jet frames}
(Definition~\ref{jet}), which record not just the sizes of the angles,
but their directions as well.

The totality of the (finite amount of) radius and angle sequence data
induces a partial order on events.  The resulting {\em source poset}\/
(Definition~\ref{poset}), which owes its existence to the finiteness
result in Theorem~\ref{frames}, describes precisely which events occur
before others---both macroscopically and infinitesimally.  Since
wavefront bifurcation is a local phenomenon at an event point,
incomparable events can occur simultaneously, or can be viewed as
occurring in any desired order.  Thus as time progresses, wavefront
expansion builds the source poset by adding one event at~a~time.

\subsection*{The algorithm}

It is one thing to order the set of events after having been given
all of them; but it is quite another to predict the ``next'' event
having been given only past events.  That the appropriate event to
add can be detected locally, and {\em without knowing future
events}, is the content of Theorem~\ref{src}.  Its importance is
augmented by it being the essential tool in making our algorithm
for constructing the source poset, and hence also the source
unfolding (Algorithm~\ref{code}).  Surprisingly, our geometric
analysis of infinitesimal wavefront expansion in
Sections~\ref{wavefront} and~\ref{construct} allows us to remove
all calculus from Theorem~\ref{src} and hence Algorithm~\ref{code}:
detecting the next event requires only standard tools from linear
algebra.

As we mentioned earlier, our original motivation for this paper was
its algorithmic applications.  Using the theoretical definitions
and results in earlier sections, we present pseudocode (a
semi-formal description) for our procedure constructing source
unfoldings in Algorithm~\ref{code}.  That our algorithm provides an
{\em efficient}\/ method to compute source unfoldings is formalized
in Theorem~\ref{timing}.

There are several reasons in favor of presenting pseudocode.
First, it underscores the explicit effective nature of our
combinatorial description of the source poset in Theorem~\ref{src}.
Second, it emphasizes the simplicity of the algorithm that results
from the apparenly complicated analysis in
Sections~\ref{geodesics}--\ref{construct}; in particular, the
reader interested only in the computational aspects of this paper
can start with Section~\ref{algorithm} and proceed backwards to
read only those earlier parts of the paper addressed in the
algorithm.  Finally, the pseudocode makes Algorithm~\ref{code}
amenable to actual implementation, which would be of interest but
lies outside the scope of this~work.

\subsection*{A note on the exposition}

Proofs of statements that may seem obvious based on intuition drawn
from polyhedral sur\-faces, or even solids of dimension~$3$, demand
surprising precision in the general case.  Occasionally, the
required adjustments in definitions and lemmas, and even in
statements of theorems, were borne out only after considering
configurations in dimension~$5$ or more.  The definition of source
image is a prime example, about which we remark in
Section~\ref{limits}, in the course of analyzing where various
hypotheses (convexity, pseudomanifold, and so on) become essential.
Fortunately, once the appropriate notions have been properly
identified, the subtlety seemingly disappears: the definitions
become transparent, and the proofs remain intuitive in low
dimension.  Whenever possible, we use figures to clarify
the~exposition.

\unit{Geodesics in polyhedral boundaries}\label{geodesics}

\noindent
In this paper, a {\em convex polyhedron}\/~$F$ of dimension~$d$ is
a finite intersection of closed half-spaces in some Euclidean
space~$\rr^d$, such that~$F$ that does not lie in a proper affine
subspace of~$\rr^d$.  The polyhedron~$F$ need not be bounded, and
comes with an induced Euclidean metric.  Gluing a finite collection
of convex polyhedra by given isometries on pairs of codimension~$1$
faces yields a {\em (finite) polyhedral cell complex}\/~$S$.  More
precisely, $S$~is a regular cell complex endowed with a metric that
is piecewise Euclidean, in which every face (closed cell) is
isometric to a convex polyhedron.

The case of primary interest is when the polyhedral cell complex~$S$
equals the boundary~$\partial P$ of a convex polyhedron~$P$ of
dimension~\mbox{$d+1$} in~$\rr^{d+1}$.
\begin{conv} \label{conv} \rm
We assume that $S = \partial P$ is a polyhedral boundary in all
theorems, proofs, and algorithms from here through
Section~\ref{algorithm}.
\end{conv}
\noindent
We do not require~$P$ to be bounded, though the reader interested in
polytopes will lose very little of the flavor by restricting to that
case.  Moreover, with the exception of Lemma~\ref{short},
Proposition~\ref{warp-cut}, Corollary~\ref{c:exp}, and
Theorem~\ref{exp}, the statements of all results from here through
Section~\ref{algorithm} are worded to hold verbatim for the more
abstract class of {\em convex polyhedral pseudomanifolds}, as we shall
see in Section~\ref{pseudomanifolds}.

Denote by $\mu$ the metric on~$S$, so $\mu(a,b)$ denotes the distance
between points $a,b \in S$.  A path~$\ga \subset S$ with endpoints $a$
and~$b$ is a {\em shortest path}\/ if its length equals $\mu(a,b)$.
Since we assume~$S$ has finitely many {\em facets}\/ (maximal faces),
such length-minimizing paths exist, and are piecewise linear.
A path $\eta \subset S$ is a {\em geodesic}\/ if $\eta$ is locally a
shortest path; i.e.\ for every $z \in \eta$ that is not an endpoint
of~$\eta$, there exist points $a,b \in \eta \minus \{z\}$ such that $z
\in \ga \ssu \eta$ for some shortest path $\ga$ connecting $a$ to~$b$.

Henceforth, as~$S$ has dimension~$d$, a face of dimension~\mbox{$d-1$}
will be called a {\em ridge}.  For convenience, we say that a
point~$x$ is {\em warped}\/ if $x$ lies in the union~$S_{d-2}$ of all
faces in~$S$ of dimension at most $d-2$, and call~$x$ {\em flat}\/
otherwise.  Every flat point has a neigh\-borhood isometric to an open
subset of~$\rr^d$.

\begin{prop} \label{warp}
If a shortest path~$\ga$ in~$S$ connects two points not lying on a
common facet, then $\ga$ has no warped points in its relative
interior.
\end{prop}
\begin{proof}
For any point $w$ lying in the relative interior of~$\ga$, the
intersection of $\ga$ with some neighborhood of~$w$ consists of two
line segments $\eta$ and~$\eta'$ that are each straight with one
endpoint at~$w$, when viewed as paths in $\rr^{d+1}$.  This is a
consequence of local length-minimization and the fact that each facet
of~$P$ is isometric to a polytope in~$\rr^d$.  Moreover, if $w$
happens to lie on a ridge while $\eta$ intersects the relative
interior of some facet containing~$w$, then local length-minimization
implies that $\eta'$ is not contained in the facet containing~$\eta$.
Lemma~\ref{short} shows that $w$ does not lie in~$S_{d-2}$, so the
point $w$ is not warped.%
\end{proof}

\begin{lemma} \label{short}
Let $\eta,\eta' \subset S$ be two paths that (i)~are straight
in~$\rr^{d+1}$, (ii)~share a common warped endpoint~$w \in S_{d-2}$,
and (iii)~do not both lie in a single facet.  There exists a
neighborhood $\co$ of $w$ in~$S$ such that for every $a \in \eta \cap
\co$ and $b \in \eta' \cap \co$, the path $\eta_{ab}$ from~$a$ to~$w$
to~$b$ along $\eta$ and $\eta'$ is not a shortest path in~$S$
between~$a$ and~$b$.
\end{lemma}
\begin{proof}
Translate $P$ so that $w$ equals the origin $\0 \in \rr^{d+1}$, and
let $Q$ be the unique minimal face of~$P$ that contains~$w$.  Since
$\eta$ and $\eta'$ do not lie in a single facet, the $2$-plane $E$
spanned by $\eta$ and~$\eta'$ meets $Q$ at exactly one point,
namely~$\0$.  Since $\dim(Q) \leq d-2$, the span of $Q$ and $E$ has
dimension at most~$d$.  Choose a line $L$ whose direction is linearly
independent from the span of $Q$ and~$E$.  Then the $3$-plane $H = L +
E$ intersects $Q$ only at~$\0$.  Replacing $P$ by $P \cap H$, we can
assume that $\dim(P) = 3$, so that $d=2$; note that $\0$ is a vertex
of $H \cap P$ by construction.

Although the case $d=2$ was proved in \cite[Theorem~4.3.5,
p.~127]{Al48} (see also \cite[Lemma~4.1]{SS}), we provide a simple
argument here, for completeness.  Let $\co \subset S$ be the
neighborhood of~$w$ consisting of all points at some fixed small
distance from the vertex~$w$.  Then $\co$ can be laid flat on the
plane $\rr^2$ by slicing along~$\eta$. One of the two points in this
unfolding that glue to~$a \in \co$ connects by a straight segment in
the unfolding to the unique point corresponding to~$b$.  This straight
segment shortcuts~$\eta_{ab}$ after gluing back to~$S$.%
\end{proof}

An illustration of Lemma~\ref{short} and its proof is given in
Figure~\ref{f:pyramid} below.
\begin{figure}[hbt]
\begin{center}
\psfrag{w}{$w$}
\psfrag{e1}{$\eta$}
\psfrag{e2}{$\eta'$}
\psfrag{a}{$a$}
\psfrag{b}{$b$}
\epsfig{file=\LOCAL/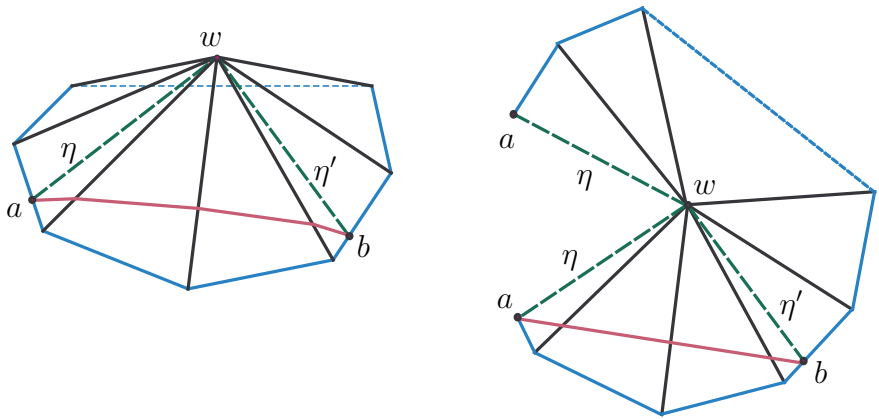,width=12cm}
\end{center}
\caption{Neighborhood of a vertex and its foldout after slicing
along the segment~$\eta$.  The points $a$ and~$b$ are connected by
a shortest~path.}
\label{f:pyramid}
\end{figure}

\begin{cor} \label{onepoint}
Let~$\eta$ be a bounded geodesic in~$S$ starting at a point~$z$ not on
any ridge.  Then $\eta$ intersects each ridge in a discrete set, so
$\eta$ traverses (in order) the interiors of a well-defined sequence
$\cl_\eta$ of facets (the\/ {\em facet sequence} of~$\eta$).
\end{cor}
\begin{proof}
Since $\eta$ is locally length minimizing, Proposition~\ref{warp}
implies that every intersection of~$\eta$ with a ridge takes place
at a flat point.  Such points have neighborhoods isometric to open
subsets of~$\rr^d$, and these intersect~$\eta$ in paths isometric
to straight segments.  It follows that $\eta$ intersects every
ridge transversely.%
\end{proof}

For each facet $F$ of~$S = \partial P$, let $\TF$ be the affine span
of~$F$ in~$\rr^{d+1}$.

\begin{defn} \label{fold}
Suppose two facets $F$ and~$F'$ share a ridge $R = F \cap F'$.  The
{\em folding map}\/ $\Phi_{F,F'} : \TF \to \TFp$ is the isometry that
identifies the copy of~$R$ in~$\TF$ with the one in~$\TFp$ in such a
way that the image of~$F$ does not intersect the interior of~$F'$.
\end{defn}

In other words, the folding map $\Phi_{F,F'}$ is the rotation of~$\TF$
with $(d-1)$-dimensional axis $R = F \cap F'$ so that $F$ becomes
coplanar with~$F'$ and lies on the other side of~$R$ from~$F'$.  It
can be convenient to view $\Phi_{F,F'}$ as rotating all of $\rr^{d+1}$
instead of only rotating $\TF$ onto~$\TFp$.  Informally, we say
$\Phi_{F,F'}$ {\em folds $\TF$ along $R$} to lie in the same affine
hyperplane as~$F'$.

\begin{defn} \label{sequential}
Given an ordered list $\cl = (F_1,F_2,\ldots,F_\ell)$ of facets such
that $F_i$ shares a (unique) ridge with $F_{i+1}$ whenever $1 \leq i <
\ell$, we write
\begin{eqnarray*}
  \Phi^{-1}_\cl &=& \Phi_{F_1,F_2}^{-1} \circ \Phi_{F_2,F_3}^{-1} \circ
  \cdots \circ \Phi_{F_{\ell-1},F_\ell}^{-1}
\end{eqnarray*}
for the {\em unfolding of $\TFn \ell$ onto~$\TFn 1$}, noting that indeed
$\Phi^{-1}_\cl(\TFn \ell) = \TFn 1$.  Setting $\cl_i =
(F_1,\ldots,F_i)$, the {\em sequential unfolding}\/ of a subset $\Ga
\subseteq F_1 \cup \cdots \cup F_\ell$ {\em along $\cl$} is the set
\begin{eqnarray*}
  (\Ga \cap F_1) \cup \Phi_{\cl_2}^{-1}(\Ga \cap F_2) \cup \cdots \cup
  \Phi_{\cl_\ell}^{-1}(\Ga \cap F_\ell) &\subset& \TFn 1,
\end{eqnarray*}
\end{defn}

By Corollary~\ref{onepoint}, we can sequentially unfold any geodesic.
Next, we use this unfolding to show uniqueness of shortest paths
traversing given facet sequences.

\begin{lemma} \label{homotopic}
Let $v$ and $w$ be flat points in~$S$.  Given a sequence $\cl$ of
facets, there can be at most one shortest path $\ga$ connecting $v$
to~$w$ such that $\ga$ traverses~$\cl_\ga = \cl$.
\end{lemma}
\begin{proof}
Let $\ga$ be a shortest path from $v$ to~$w$ traversing~$\cl$.  Inside
the union of facets appearing in~$\cl$, the relative interior of~$\ga$
has a neighborhood isometric to an open subset of~$\rr^d$ by
Proposition~\ref{warp} and the fact that the set of warped points is
closed.  Sequential unfolding of~$\ga$ into~$\TF$ for the first
facet~$F$ in~$\cl$ thus yields a straight segment in~$\TF$.  This
identifies~$\ga$ uniquely as the path in~$S$ whose sequential folding
along~$\cl$ is the straight segment in~$\TF$ connecting $v$
to~$\Phi_\cl^{-1}(w) \in \TF$.%
\end{proof}

In the proof of Lemma~\ref{homotopic}, we do not claim that the union
of facets in the list~$\cl$ unfolds sequentially without overlapping,
even though some shortest path~$\ga$ traverses~$\cl$.  However, some
neighborhood of~$\ga$ in this union of facets unfolds without
overlapping.

\begin{example}\label{brick}
Consider the unfolding of a $1 \times 1 \times 3$ rectangular box
as in Figure~\ref{f:brick}.  Denote by $F_\bot, F_\top,F_\fro,
F_\back,F_\lef,F_\rig$ the bottom, top, front, back, left, and
right faces, respectively.  Denote by~$\cl_i$ the list of facets
along which the points marked by~$i$ have been sequentially
unfolded to create the foldout~$U \ssu T_{F_\bot}$ in
Figure~\ref{f:brick}.  Then:
$$
\begin{array}{@{}rcl@{\quad}rcl@{\quad}rcl@{}}
\cl_1 &=& (F_\bot),&
\cl_2 &=& (F_\bot,F_\back),&
\cl_3 &=&(F_\bot,F_\fro),
\\
\cl_4 &=& (F_\bot,F_\back,F_\top),&
\cl_5 &=& (F_\bot,F_\fro,F_\top),&
\cl_6 &=& (F_\bot,F_\back,F_\lef),
\\
\cl_7 &=& (F_\bot,F_\fro,F_\lef),&
\cl_8 &=& (F_\bot,F_\back,F_\rig),&
\cl_9 &=& (F_\bot,F_\fro,F_\rig).
\end{array}
$$
\begin{figure}[hbt]
\begin{center}
\epsfig{file=\LOCAL/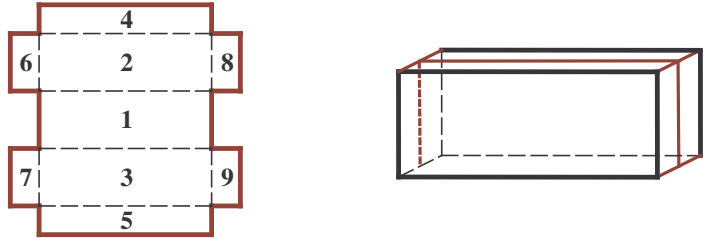,width=9.6cm}
\end{center}
\caption{An unfolding of a $1 \times 1 \times 3$ box.}
\label{f:brick}
\end{figure}
\end{example}

\unit{Cut loci}\label{cutloci}

\noindent
Most of this paper concerns the set of shortest paths with one
endpoint fixed.

\begin{defn} \label{locus}
Fix a {\em source point}\/~$v \in S$ lying interior to some facet.
A~point $x \in S$ is a {\em cut point}\/%
    \footnote{Our usage of the term `cut locus' is standard in
    differential geometry, just as our usage of `ridge' is standard in
    polyhedral geometry.  However, these usages do not agree with
    terminology in computer science, such as in~\cite{SS,AO}: their
    `ridge points' are what we call `cut points'.  Furthermore, `cut
    points' in \cite{AO} are what we would call `points on shortest
    paths to warped~points' (when $d=2$).}
if $x$ has more than one shortest path to~$v$.  Denote the set of cut
points by~$\Kv$, and call its closure the {\em cut
locus}\/~\mbox{$\KKv \subset S$}.
\end{defn}

Here is a consequence of Proposition~\ref{warp}.

\begin{cor} \label{cut}
No shortest path in~$S$ to the source point~$v$ has a cut point in its
relative interior.
\end{cor}
\begin{proof}
Suppose $c$ is a cut point in the relative interior of a shortest path
from~$v$ to~$w$.  Replacing the path from $v$ to~$c$ with another
shortest path from $v$ to~$c$ yields a new shortest path from $v$
to~$w$.  These two paths to~$w$ meet at the flat point \mbox{$c \in
S$} by Proposition~\ref{warp}.  The resulting Y-shaped intersection
at~$c$ can be improved upon in a neighborhood of $c$ isometric to an
open set in~$\rr^d$
(Fig.~\ref{f:Y}),~a~contradiction.%
\begin{figure}[hbt]
\begin{center}
\psfrag{g}{$\ga$}
\psfrag{r}{$c$}
\psfrag{V}{to~$v$}
\psfrag{W}{to~$w$}
\epsfig{file=\LOCAL/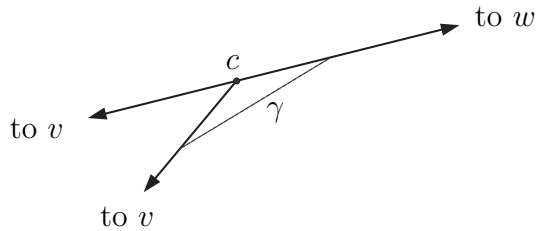,width=7cm}
\end{center}
\caption{An intersection that is Y-shaped cannot locally minimize
length in~$\rr^d$ (segment $\ga$ is a shortcut).}
\label{f:Y}
\end{figure}
\end{proof}

Our study of polyhedrality of cut loci will use Voronoi diagrams
applied to sets of points from the forthcoming definition, around
which the rest of the paper revolves.

\begin{defn} \label{d:src}
Suppose that the source point~$v$ connects by a shortest path~$\ga$
to a point~$x$ that lies on a facet~$F$ or on one of its ridges $R
\subset F$, but not on any face of~$S$ of dimension~$d-2$ or less.
If~the sequential unfolding of~$\ga$ into~$\TF$ is the segment
$[\nu,x]$, then $\nu \in \TF$ is called a {\em source image}\/
for~$F$.  Let $\src_F$ be the set of \mbox{source images for~$F$}.
\end{defn}

\begin{lemma} \label{finite}
The set\/ $\src_F$ of source images for any facet~$F$ of\/~$S$ is
finite.
\end{lemma}
\begin{proof}
The shortest path in~$\rr^{d+1}$ between any pair of distinct points
$x$ and~$y$ in a facet~$F$ is the straight segment $[x,y]$.  Since
this segment is actually contained in~$S$, any shortest path~$\ga$
in~$S$ must contain $[x,y]$ whenever it contains both~$x$ and~$y$.
Taking $x$ and~$y$ to be the first and last points of intersection
between $\ga$ and the facet~$F$, we find that~$F$ can appear at most
once in the facet sequence of a shortest path starting at the source
point~$v$.  Hence there are only finitely many possible facet
sequences of shortest paths in~$S$.  Now apply Lemma~\ref{homotopic}.%
\end{proof}

\begin{example} \label{cube-top}
Consider a unit cube with a source point in its bottom face, as in
Fig~\ref{f:cube-top}.  Then the top face has 12 source images,
shown in Fig~\ref{f:cube-top}.  The four stars~``$\star$'' are
sequential unfoldings of the source point (along three ridges each)
that are not source images: each point in the top face is closer to
some source image than to any of these stars.
\begin{figure}[hbt]
\begin{center}
\psfrag{v}{$v$}
\epsfig{file=\LOCAL/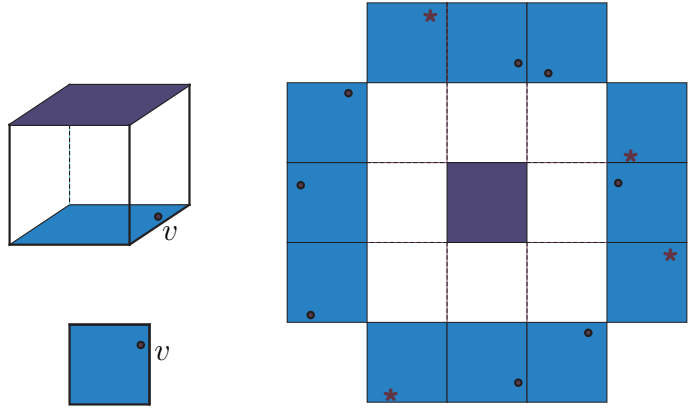,width=9.3cm}
\end{center}
\caption{Source point~$v$ on the `bottom' face, 12 source images
for the `top' face of a cube and 4 `false' source images (view from the top). }
\label{f:cube-top}
\end{figure}
\end{example}

The next result on the way to Theorem~\ref{vor} generalizes
\cite[Lemma~3.1]{Mount85} to arbitrary dimension.  Its proof is
complicated somewhat by the fact (overlooked in the proof of
\cite[Lemma~3.1]{Mount85}\footnote{
Much of \cite{Mount85}, but not Lemma~3.1, was later incorporated
and published in~\cite{MMP}.
})
that straight segments can lie inside the
cut locus, and our lack of {\sl a~priori}\/ knowledge that the cut
locus is polyhedral.

\begin{prop}[Generalized Mount's lemma] \label{mount}
Let~$F$ be a facet of~$S$, and suppose that $\nu \in \src_F$ is a
source image.  If\/ $w \in F$, then the straight segment\/ $[\nu,w]
\subset \TF$ has length at least~$\mu(v,w)$.
\end{prop}

\begin{example}
In
\begin{figure}[hbt]
\begin{center}
\psfrag{v}{$\nu$}
\psfrag{v1}{$\nu'$}
\psfrag{v2}{$\nu''$}
\psfrag{v3}{$\nu'''$}

\epsfig{file=\LOCAL/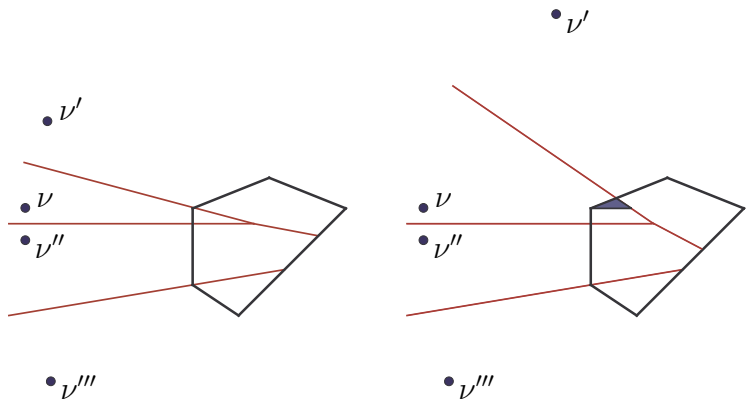,width=10cm}
\end{center}
\caption{Generalized Mount's lemma (fails for the shaded region)}
\label{f:mount}
\end{figure}
Fig.~\ref{f:mount}, the left figure is a typical illustration of
Proposition~\ref{mount} in dimension $d=2$.  In contrast, the right
figure will never occur: any point $w$ interior to the shaded
region is closest to the source image~$\nu$, but the straight
segment connecting $w$ to~$\nu$
has not been sequentially unfolded along the correct facet
sequence.
\end{example}

\begin{proof}
Since the two functions $F \to \rr$ sending $w$ to $\mu(v,w)$ and to
the length of~$[\nu,w]$ are continuous, we can restrict our attention
to those points~$w$ lying in any dense subset of~$F$.  In particular,
the cut locus has dense complement in~$F$ (Corollary~\ref{cut}) as
does the boundary of~$F$, so we assume throughout that $w$ lies in
neither the cut locus nor the boundary of~$F$.

Having fixed $\nu \in \src_F$, choose a point $x \in F$ as in
Definition~\ref{d:src}.  The set $\src_F([x,w])$ of source images
sequentially unfolded from shortest paths that end inside the segment
$[x,w]$ is finite by Lemma~\ref{finite}.  Hence we may furthermore
assume that $w$ does not lie on any hyperplane $H$ that is equidistant
from $\nu$ and a source image $\nu' \in \src_F([x,w])$.
In other words, we assume $w$ does not lie inside the hyperplane
perpendicularly bisecting any segment~$[\nu,\nu']$.

\begin{claim} \label{closer}
With these hypotheses, if $\nu \in \src_F$ but no shortest path
unfolds sequentially to the segment $[\nu,w]$, then $w$ is closer to
some point $\nu' \in \src_F([x,w])$ than~to~$\nu$.
\end{claim}

Assuming this claim for the moment, we may replace $\nu$ with~$\nu'$
and $x$ with another point~$x'$ on~$[x,w]$.  Repeating this process
and again using that the set of source images sequentially unfolded
from shortest paths ending in $[x,w]$ is finite, we eventually find
that the unique source image $\om \in \src_F([x,w])$ closest to~$w$ is
closer to~$w$ than~$\nu$~is.  Since $[\om,w]$ has length $\mu(\om,w)$,
it suffices to prove Claim~\ref{closer}.

Consider the straight segment $[x,w]$, which is contained in~$F$ by
convexity.  Let~$Y$ be the set of points $y \in [x,w]$ having a
shortest path $\ga_y$ from~$v$ that sequentially unfolds to a segment
in~$\TF$ with endpoint~$\nu$.  Then $Y$ is closed because any limit of
shortest paths from $v$ traversing a fixed facet sequence~$\cl$ is a
shortest path that sequentially unfolds along~$\cl$ to a straight
segment from the corresponding source image.
Thus, going from~$x$ to~$w$, there is a last point $x' \in Y$.  This
point~$x'$ is by assumption not equal to~$w$, so $x'$ must be a cut
point (possibly $x = x'$).

There is a facet sequence $\cl$ and a neighborhood $\co$ of~$x'$ in
$[x',w]$ such that every point in~$\co$ connects to~$v$ by a shortest
path traversing~$\cl$,
and such that unfolding the source along~$\cl$ yields a source image
$\nu' \neq \nu$ in~$\TF$.  This point~$\nu'$ connects to~$x'$ by a
segment of length $\mu(v,x')$, so the hyperplane $H$ perpendicularly
bisecting $[\nu,\nu']$ intersects $[x,w]$ at~$x'$.  By hypothesis $w
\not\in H$, and it remains to show that $w$ lies on the side of~$H$
closer to~$\nu'$.

The shortest path from $v$ to~$x'$ has a neighborhood in~$S$ disjoint
from the set of warped points and hence isometric to an open subset
of~$\rr^d$ by Proposition~\ref{warp}, because $x'$ is itself not a
warped point (we assumed $x$ lies interior to~$F$ or to a ridge~$R
\subset F$).  After shrinking~$\co$ if necessary, we can therefore
ensure that each segment $[\nu,y]$ for $y \in \co$ is the sequential
unfolding of a geodesic $\eta_y$ in~$S$.  The geodesic $\eta_y$ for $y
\in \co \minus x'$ cannot be a shortest path by definition of~$x'$, so
$[\nu,y]$ has length strictly greater than $\mu(v,y)$.  We conclude
that $\co \minus x'$, and hence also $w$, lies strictly closer to
$\nu'$ than to~$\nu$.  This finishes the proof of Claim~\ref{closer}
and with it Proposition~\ref{mount}.%
\end{proof}

\begin{figure}[hbt]
\begin{center}
\epsfig{file=\LOCAL/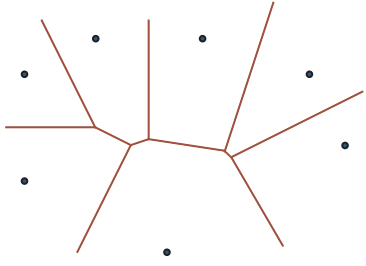,width=5cm}
\end{center}
\caption{An example of a Voronoi diagram.}
\label{f:voronoi}
\end{figure}

Before stating the first main result of the paper, we recall the
standard notion of {\em Voronoi diagram}\/ $\VV(\Ups)$ for a closed
discrete set $\Ups = \{\nu,\nu',\ldots\}$ of points in $\rr^d$.  This
is the subdivision of~$\rr^d$ whose closed cells are the sets
\begin{eqnarray*}
  V(\Ups,\nu) &=& \bigl\{\zeta \in \rr^d \:\big|\: \hbox{every point }
  \nu' \in \Ups \hbox{ satisfies } |\zeta - \nu| \leq |\zeta -
  \nu'|\bigr\}.
\end{eqnarray*}
Thus $\zeta$ lies in the interior of~$V(\Ups,\nu)$ if $\zeta$ is
closer to $\nu$ than to any other point in~$\Ups$.  An example is
given in Fig.~\ref{f:voronoi}.

\begin{thm} \label{vor}
Fix a facet~$F$ of\/~$S$, and let $V_{d-1} \subseteq \TF$ be the
union of the closed cells of dimension~\mbox{$d-1$} in the Voronoi
diagram $\VV(\src_{F\!})$ for the set of source images in~$\TF$.
If~$F^\circ$ is the relative interior of~$F$, then the set $F^\circ
\cap \Kv$ of cut points in~$F^\circ$ coincides with the
intersection $F^\circ \cap V_{d-1}$.  Moreover, if $R^\circ$ is the
relative interior of a ridge $R \subset F$, then the set $R^\circ
\cap \Kv$ of cut points in~$R^\circ$ coincides with $R^\circ \cap
V_{d-1}$.
\end{thm}

\begin{proof}
Every shortest path from the source~$v$ to a point~$w$ in $F^\circ$
or~$R^\circ$ unfolds to a straight segment in~$\TF$ of length
$\mu(v,w)$ ending at a source image for~$F$.  Proposition~\ref{mount}
therefore says that $w$ lies in the Voronoi cell $V(\Ups,\nu)$ if and
only if the segment $[\nu,w]$ has length exactly $\mu(v,w)$.  In
particular, $v$ has at least two shortest paths to~$w$ if and only if
$w$ lies in two such Voronoi cells---that is, $w \in V_{d-1}$.%
\end{proof}

To illustrate Theorem~\ref{vor} consider Example~\ref{cube-top}.
The Voronoi diagram of source images gives the cut locus in the
top face of the cube.

\begin{figure}[hbt]
\begin{center}
\psfrag{v}{$v$}
\epsfig{file=\LOCAL/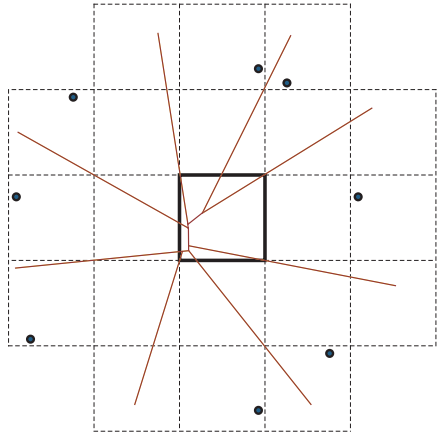,width=6.cm}
\end{center}
\caption{Cut locus of the `top' face of the cube.}
\label{f:cube-voronoi}
\end{figure}

Theorem~\ref{vor} characterizes the intersection of the cut locus with
faces of dimension $d$ or~\mbox{$d-1$} in~$S$.  For faces of smaller
dimension, we can make a blanket statement.

\begin{prop} \label{warp-cut}
Every warped point lies in the cut locus~$\KKv$; that is, $S_{d-2}
\subseteq \KKv$.
\end{prop}
\begin{proof}
It is enough to show that every point $w$ in the relative interior of
a warped face of dimension~\mbox{$d-2$} is either a cut point or a
limit of cut points, because the cut locus~$\KKv$ is closed by
definition.  Let~$\ga$ be a shortest path from~$w$ to~$v$.

First assume that every neighborhood of~$w$ contains a point having no
shortest path to~$v$ that is a deformation of~$\ga$.  Suppose that
$(y_i)_{i \in \NN}$ is a sequence of such points approaching~$w$, with
shortest paths $(\ga_i)_{i \in \NN}$ connecting the points~$y_i$
to~$v$.  Since there are only finitely many facets containing~$w$ and
finitely many source images for each facet, we may assume (by choosing
a subsequence if necessary) that for all $i$, the sequential unfolding
of~$\ga_i$ connects to the same source image for the same facet.  The
paths $(\ga_i)_{i \in \NN}$ then converge to a shortest path $\ga'
\neq \ga$ to~$w$ from~$v$, so $w \in \Kv$.

Now assume that every point in some neighborhood of~$w$ has a shortest
path to~$v$ that is a deformation of~$\ga$.  Every point on~$\ga$
other than $w$ itself is flat in~$S$ by Proposition~\ref{warp}.
Therefore some neighborhood of $\ga$ in~$S$ is isometric to an open
subset of a product $\rr^{d-2} \times C$, where $C$ is a
$2$-dimensional surface that is flat everywhere except at one point $c
\in C$ (so $C$ is the boundary of a right circular cone with
apex~$c$).  The set of points in $\rr^{d-2} \times C$ having multiple
geodesics to the image of $v$ in $\rr^{d-2} \times C$ is a relatively
open half-space of dimension $d-1$ whose boundary is $\rr^{d-2} \times
\{c\}$.  Some sequence in this open half-space converges to the image
of~$w$.
\end{proof}

Theorem~\ref{vor} and Proposition~\ref{warp-cut} imply the
following description of cut loci.

\begin{cor} \label{c:exp}
If $v$ is a source point in~$S$, then
\begin{numbered}
\item \label{item:poly}
the cut locus $\KKv$ is a polyhedral complex of dimension $d-1$, and

\item \label{item:ridge}
the cut locus $\KKv$ is the union $\Kv \cup S_{d-2}$ of the cut points
and warped points.
\end{numbered}
\end{cor}
\begin{proof}
Part~\ref{item:ridge} is a consequence of Theorem~\ref{vor} and
Proposition~\ref{warp-cut}, the latter taking care of~$S_{d-2}$, and
the former showing that points in the cut locus but outside
of~$S_{d-2}$ are in fact cut points.  Theorem~\ref{vor} also implies
that the intersection of the cut locus~$\KKv$ with any closed
facet~$F$ equals the polyhedral complex $F \cap V_{d-1}$ of
dimension~$d-1$.  Using Proposition~\ref{warp-cut} in addition to the
last sentence of Theorem~\ref{vor}, we conclude that these polyhedral
complexes glue to form a polyhedral complex of dimension~\mbox{$d-1$}
that equals the cut locus~$\KKv$, proving part~\ref{item:poly}.
\end{proof}

\unit{Polyhedral nonoverlapping unfolding}\label{unfold}

\noindent
In this section we again abide by Convention~\ref{conv}, so $S$ is a
the boundary of convex polyhedron~$P$ of dimension~\mbox{$d+1$}
in~$\rr^{d+1}$.

\begin{defn} \label{d:unfold}
A polyhedral complex $K \ssu S$ of \mbox{dimension~$d-1$} is a {\em
cut set}\/ if $K$ contains the union~$S_{d-2}$ of all closed faces of
dimension $d-2$, and $S \minus K$ is open and contractible.  A {\em
polyhedral unfolding}\/ of $S$ into $\rr^d$ is a choice of cut set $K$
and a map $S \minus K \to \rr^d$\/ that is an isometry locally on $S
\minus K$.  A {\em nonoverlapping foldout}\/ of $S$ is a surjective
piecewise linear map $\vp : \ol U \to S$ such that
\begin{numbered}
\item
$\ol U$ is the closure of its interior~$U$, which is an open
topological ball in $\rr^d$, and

\item
the restriction of $\vp$ to~$U$ is an isometry onto its image.
\end{numbered} \setcounter{separated}{\value{enumi}}
\end{defn}

Note that $K$ is not required to a polyhedral subcomplex of~$S$, but
only a subset that happens also to be a polyhedral complex; thus $K$
can `slice through interiors of facets'.  The open ball~$U$ in item~1
of the definition is usually nonconvex.  The polyhedron $P$ is a
polytope if and only if $\ol U$ is a closed ball---that is, bounded.

When the domain $\ol U$ of a nonoverlapping unfolding happens to be
polyhedral, meaning that its boundary $\ol U \minus U$ is a
polyhedral complex, the image \mbox{$K = \vp(\ol U \minus U)$} is
automatically a cut set in~$S$.  Indeed, piecewise linearity of
$\vp$ implies that $K$~is a polyhedral complex of
dimension~\mbox{$d-1$}; while the isometry implies that $K$
contains~$S_{d-2}$, and that the open ball $U \cong S \minus K$ is
contractible.  Therefore:

\begin{lemma}
If\/~$\ol U$ is polyhedral, then a nonoverlapping foldout $\vp : \ol U
\to S$ yields an ordinary polyhedral unfolding by taking the inverse
of the restriction of~$\vp$ to~$U$.
\end{lemma}

This renders unambiguous the term {\em polyhedral nonoverlapping
unfolding}.

The points in~$S$ outside of the $(d-2)$-skeleton~$S_{d-2}$ constitute
a noncompact flat Riemannian manifold~$S^\circ$.  When a point~$w$
lies relative interior to a facet~$F$, the tangent space~$\Tw$ is
identified with the tangent hyperplane~$\TF$ of~$F$, but when $w$ lies
on a ridge, there is no canonical model for~$\Tw$.

Most tangent vectors $\zeta \in \Tw$ can be exponentiated to get a
point $\exp(\zeta) \in S^\circ$ by the usual exponential map from the
tangent space~$\Tw$ to the Riemannian manifold~$S^\circ$.  (One can
show that the set of tangent vectors that cannot be exponentiated has
measure zero in~$\Tw$; we shall not use this~fact.)  In the present
case, we have a partial compactification~$S$ of~$S^\circ$, which
allows us to extend this exponential map slightly.

\begin{defn} \label{canbeexp}
Fix a point $w \in S^\circ = S \minus S_{d-2}$.  A~tangent vector
$\zeta \in \Tw$ {\em can be exponentiated}\/ if the usual exponential
of~$t\zeta$ exists in~$S^\circ$ for all real numbers~$t$ satisfying $0
\leq t < 1$.  In this case, set $\exp(\zeta) = \lim_{t \to 1}
\exp(t\zeta)$.
\end{defn}

The exponential map~\mbox{$f_\zeta: t \to \exp(t \zeta)$} takes the
interval $[0,1]$ to a geodesic~$\eta \subset S$, and should be thought
of as `geodesic flow' away from~$w$ with tangent~$\zeta$.

Henceforth fix a source point~$v \in S$ not lying on any face of
dimension less than~$d$.

\begin{defn} \label{U}
The {\em source interior}\/ $\Uv$ consists of the tangent vectors
$\zeta \in \Tv$ at the source point~$v$ that can be exponentiated, and
such that the exponentials\/ \mbox{$\exp(t \zeta)$} for\/ $0 \leq t
\leq 1$ do not lie in the cut locus~$\KKv$.  The closure of~$\Uv$ is
the {\em source foldout}\/~$\UUv$.
\end{defn}

Our next main result justifies the terminology for~$\Uv$ and its
closure~$\UUv$.

\begin{thm} \label{exp}
Fix a source point~$v$ in~$S$.  The exponential map\/ $\exp: \UUv \to
S$ from the source foldout to~$S$ is a polyhedral nonoverlapping
foldout, and the boundary $\UUv \minus \Uv$ maps onto the cut
locus~$\KKv$.  Hence $\KKv$ is a cut set inducing a polyhedral
nonoverlapping unfolding $S\minus\KKv \to \Uv$ to the source interior.
\end{thm}
\begin{proof}
It suffices to show the following, in view of parts~\ref{item:poly}
and~\ref{item:ridge} from Corollary~\ref{c:exp}.
\begin{numbered} \setcounter{enumi}{\value{separated}}
\item \label{item:star}
The metric space $S \minus \KKv$ is homeomorphic to an open ball.

\item \label{item:surj}
The exponential map $\exp: \UUv \to S$ is piecewise linear and
surjective.

\item \label{item:isom}
The exponential map $\exp: \Uv \to S \minus \KKv$ is an isometry.
\end{numbered}

Every shortest path is the exponential image of some ray in~$\UUv$ by
Proposition~\ref{warp}, and the set of vectors $\zeta \in \UUv$
mapping to $S \minus \KKv$ is star-shaped by part~\ref{item:ridge}
along with Proposition~\ref{warp} and Corollary~\ref{cut}.  This
implies part~\ref{item:star}
and surjectivity in part~\ref{item:surj}.  The space $S^\circ = S
\minus S_{d-2}$ is isometric to a flat Riemannian manifold.  Hence the
exponential map is a local isometry on any open set of tangent vectors
where it is defined.  The definition of~$\UUv$ implies that $\exp$ is
injective on the interior~$\Uv$, so the surjectivity in
part~\ref{item:surj} shows that $\exp : \Uv \to S \minus \KKv$ is an
isomorphism~of Riemannian manifolds, proving part~\ref{item:isom}.
Every isometry between two open subsets of affine spaces is linear, so
the piecewise linearity in part~\ref{item:surj} is a consequence of
part~\ref{item:isom}.%
\end{proof}

\begin{example} \label{cube-source}
Consider a cube~$P$ and a source point~$v$ located off-center on
the bottom face of~$P$, as in Example~\ref{cube-top} and
Figure~\ref{f:cube-voronoi}.  The cut locus~$\KKv$ and the
corresponding source foldout~$\UUv$ are shown in
Fig.~\ref{f:cube-source}.
See Fig.~\ref{f:two-fold} for the case
when~$v$ is in the center of the bottom face.
\begin{figure}[hbt]
\begin{center}
\psfrag{F}{$F$}
\psfrag{v}{$v$}
\epsfig{file=\LOCAL/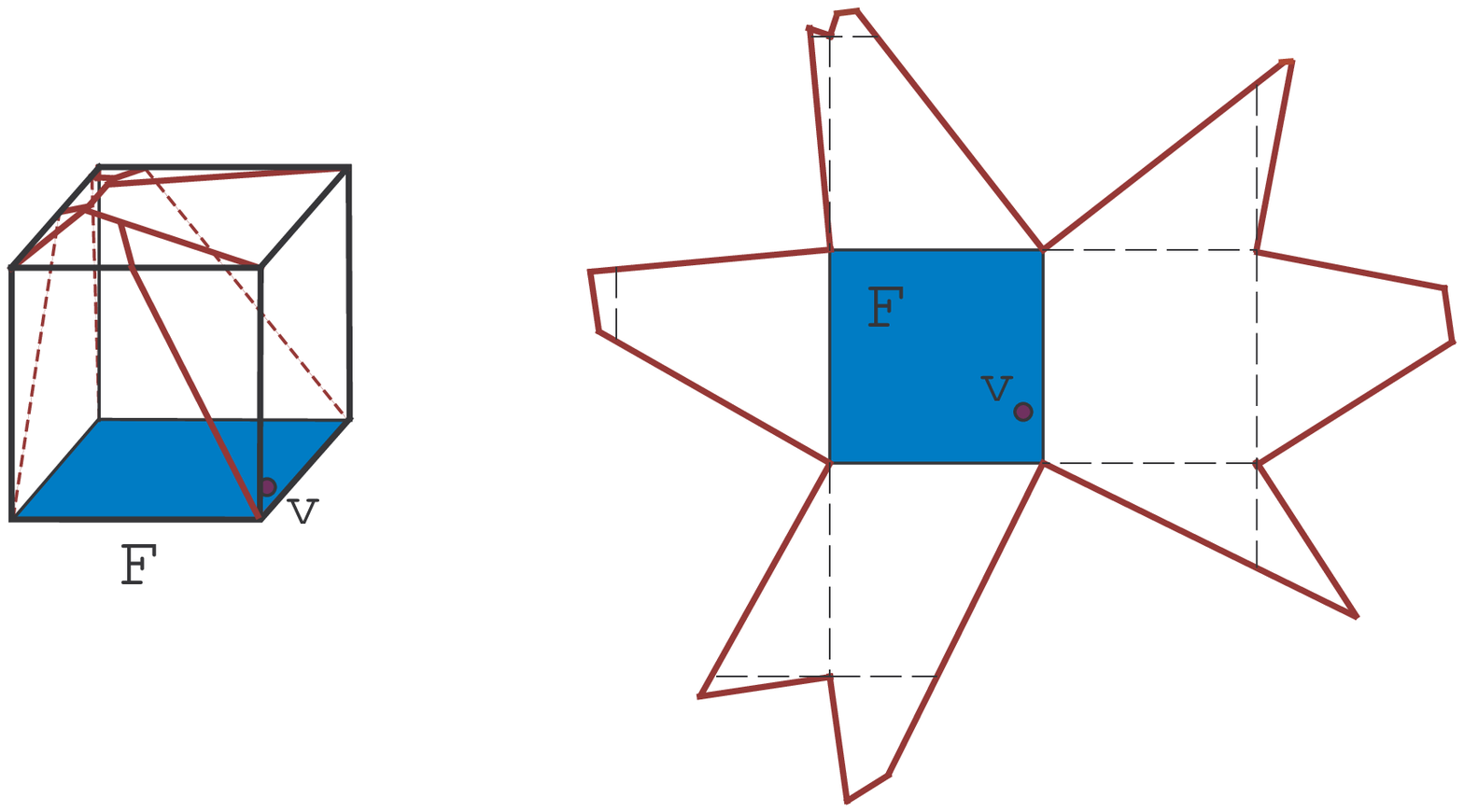,width=8.8cm}
\end{center}
\caption{Cut locus~$\KKv$ and source foldout~$\UUv$ of the
cube}
\label{f:cube-source}
\end{figure}
\end{example}

\begin{remark}
Surjectivity of the exponential map does not follow from~$S^\circ$
being a Riemannian manifold: convexity plays a crucial role (see
Fig.~\ref{non-conv}
\begin{figure}[hbt]
\begin{center}
\psfrag{v}{$v$}
\epsfig{file=\LOCAL/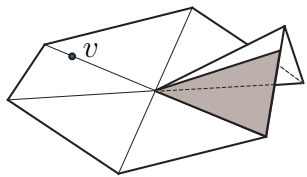,width=4.2cm}
\end{center}
\caption{Shaded region lies outside of $\exp(\Tv)$.}
\label{non-conv}
\end{figure}
for the case of a nonconvex surface.)  In fact, surjectivity
of~$\exp$ on a polyhedral manifold is equivalent---in any
dimension---to the manifold having positive curvature
\cite[Lemma~5.1]{Sto76}.  Theorem~\ref{exp} extends to the class of
{\em convex polyhedral pseudomanifolds}, but not quite verbatim;
see Theorem~\ref{pseudoexp} for the few requisite modifications.
\end{remark}

\unit{The source poset}\label{wavefront}

\noindent
In this section we define the source poset (Definition~\ref{poset}),
and in the next, we show how to build it step by step
(Theorem~\ref{src}).  The reader should consider
Definition~\ref{poset} as the main result in this section, although it
is the existence and finiteness properties for minimal jet frames in
Theorem~\ref{frames} that endow the source poset with its power to
make continuous wavefront expansion combinatorially tractable.

\begin{defn} \label{jet}
Fix a polyhedron $V$ in~$\rr^d$.  Given a list $\bar\zeta =
(\zeta_1,\ldots,\zeta_r)$ of mutually orthogonal unit vectors
in~$\rr^d$, define for $\ve \in \rr$ the unit vector
\begin{eqnarray*}
  J_{\bar\zeta}(\ve) &=& \frac{\ve\zeta_1 + \cdots + \ve^r\zeta_r}
  {\sqrt{\ve^2 + \ve^4 + \cdots + \ve^{2r}}}.
\end{eqnarray*}
If $x \in V$ and $x + \ve J_{\bar\zeta}(\ve)$ lies in~$V$ for all small
$\ve > 0$, then the vector-valued function~$J_{\bar\zeta}$\/ is a {\em
unit jet}\/ of {\em order}\/~$r$ at~$x$ in~$V$, and $\bar\zeta$ is a
{\em partial jet frame}\/ at~$x$ along~$V$.  If,~in addition, $x + \ve
J_{\bar\zeta}(\ve)$ lies relative interior to~$V$ for all small $\ve >
0$, then $\bar\zeta$ is a {\em jet frame}.
\end{defn}

The definition will be used later in the case where the convex
polyhedron~$V$ is a closed Voronoi cell $R \cap V(\src_F,\om)$ for
some ridge~$R$ of a facet~$F$, and $\om \in \src_F$ is a source image.
Think of the point $x \in V$ as the closest point in~$V$ to~$\om$.  It
will be important later (but for now may help in understanding the
next definition) to note that the relative interior of a polyhedron $V
= R \cap V(\src_F,\om)$ is contained in the relative interior of the
ridge~$R$ by Definition~\ref{d:src} and Theorem~\ref{vor}.

We do not assume the polyhedron~$V$ has dimension~$d$.  However, the
order~$r$ of a unit jet in~$V$, or equivalently the order of a jet
frame along~$V$, is bounded above by the dimension of~$V$.  In
particular, we allow $\dim(V) = 0$, in which case the only jet frame
is empty---that is, a list~$\nothing$ of length zero---and $J_\nothing
\equiv 0$.

The {\em lexicographic order}\/ on real vectors $\bar a$ and $\bar b$
of varying lengths is defined~by
\begin{eqnarray*}
  (a_1,\ldots,a_r) &<& (b_1,\ldots,b_s)
\end{eqnarray*}
if the first nonzero coordinate of $\bar a - \bar b$ is negative, where
by convention we set $a_i = 0$ for $i \geq r+1$ and $b_j = 0$ for $j
\geq s+1$.

\begin{defn} \label{minimal}
Fix a convex polyhedron~$V$ in~$\rr^d$, a point $x \in V$, and an {\em
outer support vector}\/ $\nu \in \rr^d$ for~$V$ at~$x$, meaning that
$\nu \cdot y \leq \nu \cdot x$ for all points $y \in V$.  A~jet
frame~$\bar\zeta$ at~$x$ along~$V$ is {\em minimal}\/ if the {\em
angle sequence} \mbox{$-(\nu\cdot\zeta_1,\ldots,\nu\cdot\zeta_r)$} is
lexicographically smaller than
$-(\nu\cdot\zeta'_1,\ldots,\nu\cdot\zeta'_{r'})$ for any jet
frame~$\bar\zeta'$ at~$x$ along~$V$.
\end{defn}

Again think of $V = R \cap V(\src_F,\om)$, with $\nu = \om - x$ being
the outer support vector.

In general, that $\nu$ is an outer support vector at~$x$ means
equivalently that $x$ is the closest point in~$V$ to~\mbox{$x + \nu$}.
Minimal jet frames $\bar\zeta$ can also be described more
geometrically: the angle formed by~$\nu$ and~$\zeta_1$ must be as
small as possible, and then the angle formed by~$\nu$ and~$\zeta_2$
must be as small as possible given the angle formed by~$\nu$
and~$\zeta_1$, and so~on.  It is worth bearing in mind that because
$\nu$ is an outer support vector, the angle formed by~$\nu$
and~$\zeta_1$ is at least~$\pi/2$ (that is, obtuse or right).

\begin{lemma} \label{circle}
If $\zeta$ and $\zeta'$ are vectors of equal length in\/~$\rr^d$, and
$\nu \in \rr^d$ is a vector satisfying $\nu \cdot \zeta \leq 0$ and
$\nu \cdot \zeta' \leq 0$, then $|\nu-\zeta| < |\nu-\zeta'|$ if and
only if $\nu \cdot \zeta > \nu \cdot \zeta'$.
\end{lemma}
\begin{proof}
Draw $\nu$ pointing away from the center of the circle containing
$\zeta$ and~$\zeta'$, with these vectors on the other side of the
diameter perpendicular to~$\nu$.  Then use the law of cosines: the
radii $\zeta$ and~$\zeta'$ have equal length, and $\nu$ has fixed
length; only the distances from~$\nu$ to~$\zeta$ and~$\zeta'$
change with the angles of $\zeta$ and~$\zeta'$ with~$\nu$ (see
\begin{figure}[ht]
\begin{center}

\psfrag{v}{$\nu$}
\psfrag{k}{$\zeta'$}
\psfrag{kp}{$\zeta$}
\psfrag{Q}{$|\nu - \zeta| < |\nu - \zeta'|$}

\epsfig{file=\LOCAL/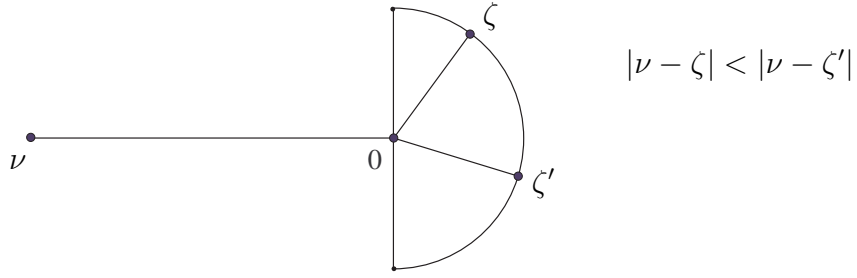,width=9cm}
\end{center}
\caption{Dot product vs.\ length}
\label{f:circle}
\end{figure}
Fig.~\ref{f:circle}).
\end{proof}

Minimal jet frames admit a useful metric characterization as follows.

\begin{prop} \label{min}
Fix two polyhedra $V$ and~$V'$ with outer support vectors $\nu$
and~$\nu'$, of equal length, at points $x \in V$ and $x' \in V'$,
respectively.  Let $\bar\zeta$ and $\bar\zeta'$ be partial jet frames
at~$x$ along~$V$ and $x'$ along~$V'$, respectively.  The angle
sequence $-\nu \cdot \bar\zeta$ is smaller than $-\nu' \cdot
\bar\zeta'$ in lexicographic order if and only if there exists $\ve_0
> 0$ such that $x+\nu$ is closer to $x+\ve J_{\bar\zeta}(\ve)$ than
$x+\nu'$ is to $x'+\ve J_{\bar\zeta'}(\ve)$ for all positive~$\ve <
\ve_0$.
\end{prop}
\begin{proof}
Since the dot product of~$\nu$ with each vector $J_{\bar\zeta}(\ve)$ or
$J_{\bar\zeta'}(\ve)$ is negative, and these are unit vectors, it is
enough by Lemma~\ref{circle} to show that minimality is equivalent~to
\begin{eqnarray*}
  \nu\cdot J_{\bar\zeta}(\ve) &\geq& \nu\cdot J_{\bar\zeta'}(\ve)
  \quad\hbox{for all nonnegative values of } \ve < \ve_0.
\end{eqnarray*}
If the first nonzero entry of $\nu\cdot\bar\zeta - \nu\cdot\bar\zeta'$
is $c = \nu\cdot(\zeta_i - \zeta_i')$, then for nonnegative values
of~$\ve$ approaching~$0$, the difference $\nu \cdot J_{\bar\zeta}(\ve)
- \nu \cdot J_{\bar\zeta'}(\ve)$ equals $c \ve^{i-1}$ times a positive
function approaching~$1$.  The desired result follows easily.%
\end{proof}

\begin{cor} \label{c:min}
Fix an outer support vector~$\nu$ at a point $x$ in a polyhedron~$V$.
A~jet frame~$\bar\zeta$ at~$x$ along~$V$ is minimal if and only if,
for every jet frame~$\bar\zeta'$ at~$x$ along~$V$, $x+\nu$ is weakly
closer to $x+\ve J_{\bar\zeta}(\ve)$ than to $x+\ve
J_{\bar\zeta'}(\ve)$ for all small nonnegative~$\ve$.
\end{cor}

It is not immediately clear from the definition that minimal jet
frames always exist: \textsl{a~priori}\/ there could be a continuum
of choices for~$\zeta_1$, and then a continuum of choices
for~$\zeta_2$ in such a way that no minimum is attained.  Although
such continua of choices can indeed occur, we shall see by
constructing minimal jet frames explicitly in Theorem~\ref{frames}
that a minimum is always attained.

First we need to know more about how (partial) jet frames at~$x$
reflect the local geometry of~$V$ near~$x$.  The {\em tangent cone}\/
to a polyhedron~$V \subseteq \rr^d$ at $x \in V$ is the~cone
\begin{eqnarray*}
  \Tx V &=& \rr_{\geq 0}\{\zeta \in \rr^d \mid x + \zeta \in V\}
\end{eqnarray*}
generated by vectors that land inside~$V$ when added to~$x$.

\begin{defn}
Fix a partial jet frame~$\bar\zeta$ at~$x$ along a polyhedron~$V$
in~$\rr^d$.  Let $\bar\zeta^\perp$ be the linear subspace of~$\rr^d$
orthogonal to the vectors in~$\bar\zeta$, and fix a sufficiently small
positive real number~$\ve$.  Then define the {\em iterated tangent
cone}\/
\begin{eqnarray*}
  \Txx{\bar\zeta} V &=& T_{\hspace{-.12ex}\xi} \big((\xi +
  \bar\zeta^\perp) \cap \Tx V\big)
\end{eqnarray*}
as the tangent cone at~$\xi = J_{\bar\zeta}(\ve)$ to the intersection of
$\Tx V$ with the affine space~\mbox{$\xi + \bar\zeta^\perp$}.
\end{defn}

Just as the partial jet frames of order~$1$ generate the tangent
cone~$\Tx V$, we have the following characterization of iterated tangent
cones.  We omit the easy proof.

\begin{lemma}
The iterated tangent cone $\Txx{\bar\zeta} V$ is generated by all
unit vectors $\zeta_{r+1}$ in\/~$\rr^d$ extending the partial jet
frame $\bar\zeta = (\zeta_1,\ldots,\zeta_r)$ to a partial jet frame
$(\zeta_1,\ldots,\zeta_r,\zeta_{r+1})$ of order~\mbox{$r+1$}.  In
particular, iterated tangent cones do not depend on the small~$\ve
> 0$.
\end{lemma}

Now we set out to construct minimal jet frames inductively.

\begin{lemma} \label{sharp}
Fix a polyhedron~$V$ and an outer support vector~$\nu$ at $x \in V$.
If~$\zeta \in \Tx V$ is a unit vector with $\nu \cdot \zeta$ maximal,
then $\nu$ is an outer support vector at\/ $\0 \in \Txx\zeta V$.
\end{lemma}
\begin{proof}
If~$\zeta' \in \Txx\zeta V$ is a unit vector satisfying $\nu \cdot
\zeta' > 0$, then $\xi = (\zeta + \ve \zeta')/\sqrt{1+\ve^2}$ for
small $\ve > 0$ is a unit vector in~$\Tx V$ satisfying $\nu \cdot \xi
> \nu \cdot \zeta$, contradicting maximality.%
\end{proof}

In `generic' cases, the functional $\zeta \mapsto \nu \cdot \zeta$ for
an outer support vector~$\nu$ on a cone will take on the maximum value
zero uniquely at the origin.  In this case, as we now show, there can
be only finitely many unit vectors~$\zeta$ in the cone having $\nu
\cdot \zeta$ maximal, and these lie along the {\em rays}, meaning
one-dimensional faces of the cone.  Note that genericity forces the
cone to be {\em sharp}, meaning that it contains no linear subspaces.

\begin{prop} \label{ray}
Let $\nu$ be an outer support vector for a sharp polyhedral cone~$C$,
and assume $\nu$ is maximized uniquely at the origin\/~$\0$.  The
minimum angle between~$\nu$ and a unit vector $\zeta \in C$ occurs
when $\zeta$ lies on a ray of~$C$.
\end{prop}
\begin{proof}
Let $Z$ be the set of unit vectors in~$C$.
Suppose that $L$ is a $2$-dimensional subspace inside the span of~$C$,
and let $\bar\nu$ be the orthogonal projection of~$\nu$ onto~$L$.  View
$\nu$ and $\bar\nu$ as functionals on $L$ via $\zeta \mapsto \nu \cdot
\zeta$, and observe that $\nu \cdot \zeta = \bar\nu \cdot \zeta$ for all
$\zeta \in L$.  The circular arc $Z \cap L$ lies inside the unit circle
in~$L$, and $\bar\nu$ takes nonpositive values on $Z \cap L$ because
$\nu$ is an outer support vector.  Elementary geometry shows that
$\bar\nu$ is therefore maximized on $Z \cap L$ only at one or both of
the endpoints of the arc~$Z \cap L$.  This argument proves that $\nu$
cannot be maximized on~$Z$ at a point $\zeta \in Z$ unless $\zeta$ lies
in the boundary of~$Z$.  The result now follows by induction on the
dimension of the cone~$C$.
\end{proof}

In `nongeneric' cases, including when the polyhedral cone~$C$ has
nonzero {\em lineality}, which is by definition the largest vector
space contained in~$C$, the functional $\nu$ is maximized along a face
of positive dimension.  In this case, there is always a continuum of
choices for unit vectors $\zeta \in C$ having $\nu \cdot \zeta = 0$.
However, the sequences of iterated tangent cones to appear in
Theorem~\ref{frames} will not in any noticeable way depend on the
continuum of choices, because of the next result.

\begin{lemma} \label{continuum}
Fix a polyhedron~$V$, a point $x \in V$, and a face $F$ of~$V$
containing~$x$.  The iterated tangent cone $\Txx{\bar\zeta} V$ is
independent of the jet frame~$\bar\zeta$ for~$F$ at~$x$.
\end{lemma}
\begin{proof}
Translate~$V$ so $x + \ve J_{\bar\zeta}(\ve)$ lies at the origin~$\0 \in
\rr^d$.  Then~$F$ spans a dimension $\dim(F)$ linear subspace $\<F\>
\subseteq \rr^d$, and the iterated tangent cone is $\Txx{\bar\zeta} V =
\<F\>^\perp \cap T_{\hspace{-.12ex}\0} V$.  Now use the fact that
$T_{\hspace{-.12ex}\0} V = \Txi V$ for all vectors $\xi$ relative
interior~to~$F$.%
\end{proof}

The main theorem in this section says that given an outer support
vector~$\nu$, there is a finite procedure using elementary linear
algebra for producing a single jet frame that is, in a precise sense,
tilted as much toward~$\nu$ as possible.

\begin{thm} \label{frames}
Fix a polyhedron~$V$ and an outer support vector~$\nu$ at $x \in V$.
Inductively construct a finite set of jet frames for~$V$ at~$x$ by
iterating the following procedure.  For each of the finitely many
partial jet frames $\bar\zeta$ already constructed:
\begin{itemize}
\item
If $\nu$ is orthogonal to a nonzero vector in~$\Txx{\bar\zeta} V$,
then add any such vector to~$\bar\zeta$.

\item
If $\nu \cdot \zeta < 0$ for all nonzero vectors~$\zeta$
in~$\Txx{\bar\zeta} V$, then create one new partial jet frame for each
of the (finitely many) rays of~$\Txx{\bar\zeta} V$ minimizing the
angle with~$\nu$, by appending to~$\bar\zeta$ the unit vector along
that ray.
\end{itemize}
At least one of the finitely many jet frames constructed in this way is
minimal.
\end{thm}
\begin{proof}
The sequences of vectors constructed by the iterated procedure are jet
frames by Lemma~\ref{sharp}.  Given an arbitrary jet frame~$\bar\xi$
for~$V$ at~$x$, it is enough to show that the angle sequence
of~$\bar\xi$ satisfies $\nu \cdot \bar\zeta \geq \nu \cdot \bar\xi$ in
lexicographic order for some constructed jet frame~$\bar\zeta$.
Indeed, then a jet frame whose angle sequence is lexicographically
minimal among the constructed ones is minimal.  Suppose that the
first~$i-1$ entries $(\xi_1,\ldots,\xi_{i-1})$ agree with a
constructed jet frame, but that $(\xi_1,\ldots,\xi_i)$ do not.

If $\nu \cdot \xi_i < 0$ then $\nu \cdot \xi_i$ is less than $\nu \cdot
\zeta_i$ for some constructed jet frame~$\bar\zeta$ agreeing
with~$\bar\xi$ through the $(i-1)^\st$ entry, by Proposition~\ref{ray}.

If,~on the other hand, $\nu \cdot \xi_i = 0$, then pick the index~$j$
maximal among those satisfying $\xi_j \neq 0$ and also $\nu \cdot \xi_i
= \cdots = \nu \cdot \xi_j = 0$.  If there is a constructed jet
frame~$\bar\zeta$ that agrees with $\bar\xi$ through the~$j^\th$ entry,
but has $\nu \cdot \xi_{j+1} < \nu \cdot \zeta_{j+1} = 0$, then already
we are done.  Therefore we can assume that the constructed jet
frame~$\bar\zeta$ agrees with~$\bar\xi$ through index~\mbox{$(i-1)$},
that $\bar\zeta$ has $\nu\cdot\zeta_i = \cdots = \nu\cdot\zeta_j = 0$,
and that either $\nu\cdot\zeta_{j+1} < 0$ or else~$\bar\zeta$ has
order~$j$.  Replacing the vectors $\xi_i,\ldots,\xi_j$ in~$\bar\xi$ with
$\zeta_i,\ldots,\zeta_j$ yields a new jet frame~$\bar\xi'$, by
Lemma~\ref{continuum} applied to the face~$F$ of the iterated tangent
cone~$\Txx{\bar\zeta'} V$ orthogonal to~$\nu$, where $\bar\zeta' =
(\zeta_1,\ldots,\zeta_{i-1})$.  Downward induction on the number of
entries of~$\bar\xi'$ shared with a constructed jet frame completes the
proof.%
\end{proof}

Our goal is to apply jets to define a poset structure on the set of
source images.  First, we need some terminology and preliminary
concepts.  The next definition is made in slightly more generality
than required for dealing only with complete sets of source images
because we shall need it for Theorem~\ref{src}.

Resume the notation from previous sections regarding the polyhedral
complex~$S$.  Recall that $\TF \cong \rr^d$ is the tangent
hyperplane to the facet~$F$.  Removing from~$\TF$ the affine span
$\TR$ of any ridge $R \subset F$ leaves two connected components
(open half-spaces).  Thus it makes sense to say that a point $\nu
\in \TF \minus \TR$ lies either on the {\em same side}\/ or on the
{\em opposite side}\/ of~$R$ as does~$F$.

\begin{defn} \label{cansee}
Fix a facet~$F$, a ridge $R \subset F$, and a finite set $\Ups \subset
\TF$.
\begin{numbered}
\item
A point $\om \in \Ups$ {\em can see $F$ through~$R$ in~$\VV(\Ups)$}\/
if $\om$ lies on the opposite side of~$R$ as $F$~does, and the closed
Voronoi cell $V(\Ups,\om)$ contains a point interior to~$R$.

\item
A point $\om \in \Ups$ {\em can see $R$ through~$F$ in~$\VV(\Ups)$}\/
if $\om$ lies on the same side of $R$ as $F$~does, and the closed
Voronoi cell $V(\Ups,\om)$ contains a point interior to~$R$.

\item
In either of the above two cases, the ridge~$R$ lies at {\em
radius~$r = r(R,\om)$}\/ from~$\om$ if $r$~equals the smallest
distance in~$\TF$ from~$\om$ to a point of~$R \cap V(\Ups,\om)$.

\item
The unique {\em closest point}\/ $\rho(R,\om)$ to~$\om$ in $R \cap
V(\Ups,\om)$ has distance~$r$~from~$\om$.

\item
The {\em outer support vector}\/ of the pair $(R,\om)$ is
$\om-\rho(R,\om)$.

\item
The {\em angle sequence}\/ $\angle(R,\om)$ is the angle sequence
$-(\om - \rho(R,\om)) \cdot \bar\zeta$ for any minimal jet
frame~$\bar\zeta$ at~$\rho(R,\om)$ \mbox{along~$R \cap V(\Ups,\om)$}.
\end{numbered}
\end{defn}

\begin{example}
Fig.~\ref{f:cansee} depicts examples of the notions from
Definition~\ref{cansee}.
\begin{figure}[ht]
\begin{center}
\psfrag{w}{$\om$}
\psfrag{w1}{$\om^\pr$}

\psfrag{O}{$0$}
\psfrag{F}{$F$}

\psfrag{R}{$R$}
\psfrag{R1}{$R'$}
\psfrag{R2}{$R''$}

\psfrag{P}{$\rho(R,\om)$}
\psfrag{P1}{$\rho(R',\om)$}
\psfrag{P2}{$\rho(R'',\om)$}

\psfrag{Q}{$\om - \rho(R,\om)$}

\epsfig{file=\LOCAL/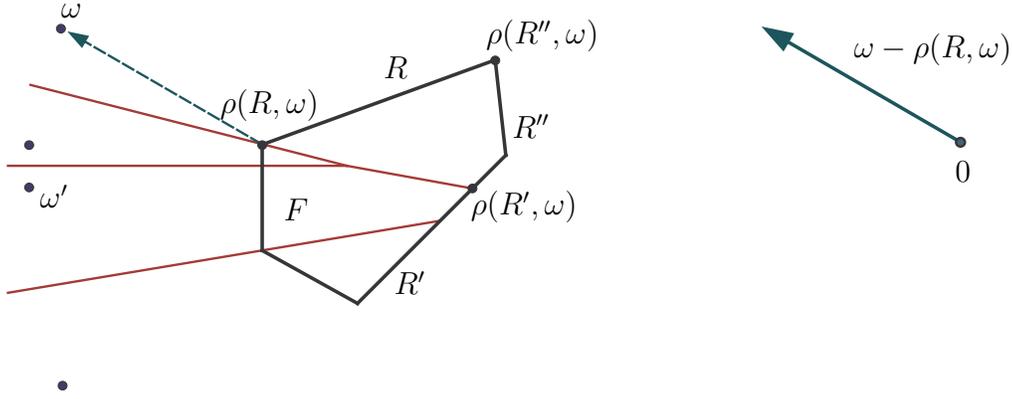,width=13.3cm}
\end{center}
\caption{Illustrations for Definition~\ref{cansee}}
\label{f:cansee}
\end{figure}
The solid pentagon is the face~$F$, while the set~$\Ups$ contains four
points.  The point $\om$ can see $F$ through the ridge~$R$, and can
see the ridges~$R'$ as well as~$R''$ through~$F$.  The three closest
points for these are indicated, as is the outer support vector
for~$(R,\om)$.  The point~$\om'$ can see the ridge~$R'$ through~$F$,
but $\om'$ can {\em not}\/ see~$R''$ through~$F$, because~$\om$ is
closer to every point of~$R''$.
\end{example}

In our applications, the finite set $\Ups$ will always be a subset of
source images in~$\src_F$, often a proper subset.  Now we are ready
for the main definition of this section.  It may help to recall that
each source image $\nu \in \src_F$ can see~$F$ through a {\em
unique}\/ ridge~$R$ by Theorem~\ref{vor}, when $\rr^d = \TF$ and the
finite set~$\Ups$ in Definition~\ref{cansee} equals~$\src_F$.

\begin{defn} \label{poset}
Fix a source point~$v$ in~$S$.  An {\em event}\/ is a pair $(\nu,F)$
with $\nu \in \src_F$ a source image for the facet~$F$.  The event
$(\nu,F)$ has
\begin{numbered}
\item
{\em radius}\/ $r(\nu,F)$ equal to the radius $r(R,\nu)$ from~$\nu$ to
the ridge~$R$ through which $\nu$ can see~$F$ in the Voronoi
subdivision $\VV(\src_{F\!})$ of~$\TF$;

\item
{\em event point}\/ $\rho(\nu,F)$ equal to the closest
point~$\rho(R,\nu)$ in~$R \cap V(\src_F,\nu)$ to~$\nu$;~and

\item
{\em angle sequence}\/~$\angle(\nu,F)$ equal to the angle sequence
$\angle(R,\nu)$.
\end{numbered}
(The trivial event $(v,\hbox{facet}(v))$ has radius~$0$, event
point~$v$, and empty angle sequence.)  The {\em source poset}\/
$\src(v,S)$ is the set of events, partially ordered with $(\nu,F) \prec
(\nu',F')$~if
\begin{itemize}
\item
$r(\nu,F) < r(\nu',F')$, or if

\item
$r(\nu,F) = r(\nu',F')$ and $\angle(\nu,F)$ is lexicographically smaller
than~$\angle(\nu',F')$.
\end{itemize}
\end{defn}

\begin{remark}
Corollary~\ref{c:min} says that breaking ties by lexicographically
comparing angle sequences at event points is the same as breaking ties
by comparing distances from each source image to a minimal jet at its
event point.  This is the precise sense in which the source poset
orders events by comparing infinitesimal expansion of the wavefront
along the interiors of ridges containing event points.
\end{remark}

\unit{Constructing source images}\label{construct}

\noindent
Aside from its abstract dynamical interpretation, the importance of
the source poset here stems from its ability to be computed
algorithmically, as we shall see here and in
Section~\ref{algorithm}.  Source images are built one by one, using
only previously built source images as stepping stones.  These
stepping stones form an {\em order ideal}\/ in $\src(v,S)$, meaning
a subset $\II \subset \src(v,S)$ closed under going down: $E \in
\II$ and \mbox{$E' \prec E$ $\implies$ $E' \in \II$}.

To make a precise statement in the main result, Theorem~\ref{src}, we
need one more dose of terminology, describing constructions in~$S$
determined by a choice of order~ideal.

\begin{defn} \label{ideal}
Fix an order ideal~$\II$ in the source poset\/~$\src(v,S)$.  For
each facet~$F$, let $\UF \subset \TF$ be the set of source images
$\om \in \src_F$ with $(\om,F) \in \II$.  The set $\EE_\II$ of {\em
potential events}\/ consists of triples $(\om,F,R')$ such that
\begin{itemize}
\item
$\om$ can see the ridge~$R'$ through~$F$ in the Voronoi diagram
$\VV(\UF)$, but

\item
a second facet~$F'$ contains~$R'$, and the unfolding $\om' =
\Phi_{F,F'}(\om)$ of~$\om$ onto the tangent space~$\TFp$ results in a
pair $(\om',F')$ that does not lie in~$\II$.
\end{itemize}
If~$(\om',F')$ is an event in $\src(v,S) \minus \II$, then we say it is
obtained by {\em processing}\/ $(\om,F,R')$.  A potential event $E \in
\EE_\II$ is {\em minimal}\/ if it has minimal radius~$r$ among potential
events, and lexicographically minimal angle sequence among potential
events with~radius~$r$.
\end{defn}

Tracing back through notation, if $E = (\nu,F,R')$ is a minimal
potential event, then the minimal radius is $r = r(R',\nu)$, and the
minimal angle sequence is~$\angle(R',\nu)$.

\begin{thm} \label{src}
Given a nonempty order ideal~$\II$ in the source poset~$\src(v,S)$,
pick a minimal potential event $(\nu,F,R')$ in~$\EE_\II$.  If~$\nu' =
\Phi_{F,F'}(\nu)$ is the unfolding of~$\nu$ to the other facet~$F'$
containing~$R'$, then $\II' = \II \cup\{(\nu',F')\}$ is an order ideal
in~$\src(v,S)$.
\end{thm}

The statement has two parts, really: first, $\nu' \in \TFp$ is indeed
a source image; and second, $\II'$~is an order ideal in the
poset~$\src(v,S)$.  To~prove the theorem we need a number of
preliminaries.  We state results requiring an order ideal inside the
source poset~$\src(v,S)$ using language that assumes an order
ideal~$\II$ has been fixed.

Recall from Section~\ref{geodesics} the notion of facet sequence
$\cl_\ga$ for a shortest path~$\ga$.  If,~on the way to a facet~$F'$,
a~shortest path~$\ga$ from the source point~$v$ traverses a facet~$F$,
then the corresponding source images in~$F$ and~$F'$ have a special
relationship.  Precisely:

\begin{defn} \label{geoprec}
Let $(\nu,F) \prec (\nu',F')$ be events in the source poset.  Suppose
some shortest path~$\ga$ has facet sequence $\cl_\ga =
(F_1,\ldots,F_{\ell'})$ with a consecutive subsequence
\begin{eqnarray*}
  \cl = (F_\ell,\ldots,F_{\ell'}) &\hbox{in which}& F = F_\ell \hbox{
  and } F' = F_{\ell'}.
\end{eqnarray*}
If $\nu' = \Phi_\cl(\nu) = \Phi_{\cl_\ga}(v)$ is the sequential
unfolding of the source along~$\ga$, and also the sequential unfolding
of $\nu \in \TF$ into~$\TFp$, then $(\nu,F)$ {\em geodesically
precedes}\/ $(\nu',F')$.  We also say that the shortest path~$\ga$
described above is {\em geodesically preceded}\/ by~$(\nu,F)$.
\end{defn}

Since the Voronoi cells in Theorem~\ref{vor} come up so often, it will
be convenient to have easy terminology and notation for them.

\begin{defn} \label{cutcell}
Given a source image $\om \in \src_F$, the {\em cut cell}\/ of~$\om$
is $V_\om = V(\src_F,\om)$.
\end{defn}

Roughly speaking, our next result says that angle sequences increase at
successive events along shortest paths, when the event point is pinned
at a fixed point~$x$.

\begin{prop} \label{jetsincrease}
If $(\nu,F)$ geodesically precedes $(\nu',F')$ then $(\nu,F) \prec
(\nu',F')$.
\end{prop}
\begin{proof}
Because of the way partial order on~$\src(v,S)$ is defined, we may as
well assume that $F$ and~$F'$ share a ridge~$R'$, and that $\nu' =
\Phi_{F,F'}(\nu)$ is obtained by folding along this ridge.  In
addition, we may as well assume that both event points $\rho(\nu,F)$
and~$\rho(\nu',F')$ equal the same point~$x \in S$, since otherwise
$r(\nu,F) < r(\nu',F')$.  Translate to assume this point~$x$ equals
the origin~$\0$, to simplify notation.  Let $R$ be the ridge through
which $\nu$ can see~$F$, and set $V = R \cap V_\nu$ and $V' = R' \cap
V_{\nu'}$; these are the cut cells through which the source images
$\nu$ and~$\nu'$ see their corresponding~facets.

The angle geometry of~$\nu'$ relative to~$V'$ in~$\TFp$ is {\em
exactly}\/ the same as the geometry of~$\nu$ relative to~$V'$ in~$\TF$,
because $\nu'$ is obtained by rotation around an axis in~$\rr^{d+1}$
containing~$V'$.  In other words, $\nu - \nu'$ is orthogonal to~$V'$.
Therefore we need only compare the angles with~$\nu$ of jets along $V$
and~$V'$.  All jet frames will be at~$x$.

Suppose the finite sequence $(\xi_1,\xi_2,\ldots)$ is a jet frame
along~$V'$.  Noting that~$V$ and $V'$ have disjoint interiors, choose
the index~$r$ so that $\bar\xi = (\xi_1,\ldots,\xi_{r-1})$ is a
partial jet frame along~$V$, but $\bar\xi' = (\xi_1,\ldots,\xi_r)$ is
not.  It is enough to demonstrate that some partial jet frame
$(\xi_1,\ldots,\xi_{r-1},\zeta_r)$ along~$V$ has lexicographically
smaller angle sequence than~$\bar\xi'$.  Equivalently, it is enough to
produce a unit vector~$\zeta_r$ in the iterated tangent cone
$\Txx{\bar\xi} V$ satisfying $\nu \cdot \zeta_r > \nu \cdot \xi_r$.

Since $R' \cap V_{\nu'} = R' \cap V_\nu$ by Theorem~\ref{vor}, every
line segment from $\nu$ to a point in~$V'$ passes through~$V$.
Therefore, since we have translated to make $x = \0$, every segment
connecting $\nu$ to~$\Tx V'$ passes through~$\Tx V$.
This observation will become crucial below; for now, note the
resulting inequality $\dim(V) \geq \dim(V')$, which implies that the
iterated tangent cone $\Txx{\bar\xi} V$ contains nonzero vectors.  All
such vectors by definition lie in the subspace $\bar\xi^\perp$
orthogonal to the space $\<\bar\xi\>$ with basis
$\xi_1,\ldots,\xi_{r-1}$.  The same holds for~$\xi_r$, so we may
replace~$\nu$ with a vector $\om \in \bar\xi^\perp$ by adding a
\mbox{vector in~$\<\bar\xi\>$, since then}
\begin{eqnarray*}
  \om \cdot \zeta &=& \nu \cdot \zeta \quad\hbox{for all vectors}\quad
  \zeta \in \bar\xi^\perp.
\end{eqnarray*}

Fix a small positive real number~$\ve$.  The line segment
$[\nu,J_{\bar\xi'}(\ve)]$ intersects~$\Tx V$ at a point
near~$J_{\bar\xi}(\ve)$.  The image segment in~$\bar\xi^\perp$ by
orthogonal projection modulo $\<\bar\xi\>$ is $[\om,\la\xi_r]$, for
$\la = \ve^r/\sqrt{\ve^2 +\cdots+ \ve^{2r}}$.  This image segment
passes through the cone~$\Txx{\bar\xi} V$
%
at some point~$\zeta$ on its way from~$\om$ to~$\la\xi_r$.  Elementary
geometry of the triangle
\begin{figure}[ht]
\begin{center}

\psfrag{z}{\raisebox{.3ex}{$\zeta$}}
\psfrag{w}{$\om$}
\psfrag{L}{$\la\xi_r$}
\psfrag{0}{$\0$}
\epsfig{file=\LOCAL/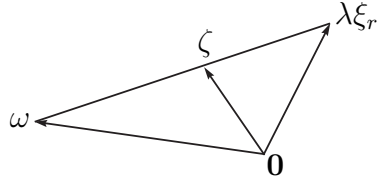,width=5cm}
\end{center}
\caption{Geodesic precedence implies smaller angle sequence}
\label{f:obtuse}
\end{figure}
with vertices $\0$, $\om$, and~$\la\xi_r$ (see Fig.~\ref{f:obtuse})
shows that the angle between $\om$ and~$\zeta$ is smaller than the
angle between $\om$ and~$\la\xi_r$.  Taking $\zeta_r = \zeta/|\zeta|$
completes the proof.%
\end{proof}

After choosing a minimal potential event~$E$, we must make sure that
when all is said and done, none of the other potential events end up
below~$E$ in the source poset.

\begin{lemma} \label{potevent}
Suppose $(\om,F,R') \in \EE_\II$ is a potential event with angle
sequence~$\angle$ and radius~$r$.  Let~$F'$ be the other facet
containing~$R'$ and $\om' = \Phi_{F,F'}(\om)$ the unfolding of $\om
\in \TF$ onto~$\TFp$.  If $(\om',F')$ is an actual event, then it
either has radius strictly bigger than~$r$, or else its angle sequence
$\angle(\om',F')$ is lexicographically larger than~$\angle$.
\end{lemma}
\begin{proof}
Assume that $(\om',F')$ is an event.  Quite simply, the result is a
consequence of the fact that the cut cell $R' \cap V_{\om'} = R' \cap
V_\om$ must be contained inside $R' \cap V(\UF,\om)$, which follows
because $\UF \subseteq \src_F$.%
\end{proof}

In comparing a newly processed event (in source poset order) to other
as yet unprocessed events, we will need to know approximately how
those other events will eventually arise.  This requires the
forthcoming lemma, in which a~{\em flat triangle}\/ inside~$S$ is any
subset of~$S$ isometric to a triangle in the Euclidean plane~$\rr^2$.

\begin{lemma} \label{Ox}
Fix a point~$x \in S$.  There is an open neighborhood\/ $\OO_x$ of~$x$
in~$S$ such that, given\/~$y \in \OO_x$ and a shortest path~$\ga$ from
the source point~$v$ to~$y$, some shortest path~$\ga'$ from~$v$ to~$x$
has the following property: the loop formed by traversing~$\ga'$ and
then the segment $[x,y]$ and finally the reverse of~$\ga$ bounds a
flat triangle in~$S$.
\end{lemma}
\begin{proof}
Choose~$\OO_x$ so small that the only closed faces of the cut
locus~$\KKv$ intersecting~$\OO_x$ are those containing~$x$.  Every cut
cell containing $y \in \OO_x$ also contains~$x$ by construction.
Convexity of cut cells (Theorem~\ref{vor}) implies that the segment
$[x,y]$ lies inside every cut cell containing~$y$ (there may be more
than one if $y$ is itself a cut point).  The source image obtained by
sequentially unfolding~$\ga$ therefore connects to every point
of~$[x,y]$ by a straight segment that sequentially folds to a shortest
path.  The union of these shortest paths is the flat triangle in
question.%
\end{proof}

Conveniently, all of the shortest paths to~$x$ already yield events
in~$\II$:

\begin{lemma} \label{closest}
Suppose some minimal potential event $E \in \EE_\II$ has closest
point~$x$.  Let $G$ be the last facet whose interior is traversed by a
shortest path~$\ga$ from the source point to~$x$.  If~$\om \in \src_G$
is the source image sequentially unfolded along~$\ga$, then $(\om,G)
\in \II$.
\end{lemma}
\begin{proof}
As $\ga$ enters~$G$, it crosses the relative interior of some ridge
of~$G$ at a point~$w$.  The event point $\rho(\om,G)$ can be no
farther than~$w$ from~$\om$.  On the other hand, $\mu(v,w) <
\mu(v,x)$, because $\ga$ traverses the interior of~$G$.  Therefore
$(\om,G)$ has radius less than $r(E) = \mu(v,x)$.%
\end{proof}


\begin{proofof}{Theorem~\ref{src}}
Suppose the minimal potential event $(\nu,F,R')$ has closest point $x
= \rho(R',\nu)$ to $R' \cap V(\UF,\nu)$, of radius~$r$, and a minimal
jet frame~$\bar\zeta$ at~$x$ with angle sequence~$\angle$.

Let~$\ga$ be a shortest path from the source that ends at a point in
the neighborhood~$\OO_x$ from Lemma~\ref{Ox}.  By that lemma and
Lemma~\ref{closest}, $\ga$~unfolds to produce a source image whose
event either lies in~$\II$, or is obtained by processing a potential
event in~$\EE_\II$, or is geodesically preceded by such a processed
event.  Applying Lemma~\ref{closest} and then
Proposition~\ref{jetsincrease}, we find that all events with event
point~$x$ that are not in~$\II$ have angle sequences lexicographically
larger than~$\angle$.

For positive~$\ve$, set $y(\ve) = x + \ve J_{\bar\zeta}(\ve)$.  When
$\ve$ is small enough, $y(\ve)$ lies interior to~$R'$, and close to~$x$,
in the neighborhood~$\OO_x$ from Lemma~\ref{Ox}.  By the previous
paragraph, every source image containing~$y$ in its cut cell is either
in~$\II$, or has angle sequence lexicographically larger than~$\angle$.

Let us now compare, for all small positive~$\ve$, the distance
to~$y(\ve)$ from~$\nu$ with the distance to~$y(\ve)$ from any source
image in $\src_F$ or~$\src_{F'}$.  Clearly the distance from a source
image~$\om$ is minimized when $y(\ve)$ lies in the cut cell~$V_\om$.
Moreover, we may restrict our attention to those source images~$\om$
whose cut cells~$V_\om$ contain~$y(\ve)$ for all sufficiently small
positive~$\ve$.
Definition~\ref{jet} says that $\bar\zeta$ is a jet frame at~$x$
along~$V_\om$.  Therefore, by Proposition~\ref{min}, we conclude using
the last sentence of the previous paragraph that $y(\ve)$ is weakly
closer to~$\nu$ than to~$\om$ for all small positive~$\ve$.  This
argument shows that $\nu'$ is a source image, so $(\nu',F')$ is an
event.  Moreover, it shows:

\begin{claim} \label{minjet}
Any minimal jet frame $\bar\zeta$ at the event point $x =
\rho(R',\nu)$ along the polyhedron $R' \cap V(\UF,\nu)$ is a minimal
jet frame at\/~$x$ along $R' \cap V_{\nu'}$.
\end{claim}

Every event in~$\src(v,S) \minus \II$ is either obtained by processing a
potential event in~$\EE_\II$, or is geodesically preceded by such a
processed potential event.  Using Claim~\ref{minjet},
we conclude by Lemma~\ref{potevent} and Proposition~\ref{jetsincrease}
that $\II'$ is an order ideal.%
\end{proofof}

\unit{Algorithm for source unfolding}\label{algorithm}

\noindent
The primary application of the analysis up to this point is an
algorithmic construction of nonoverlapping unfoldings of convex
polyhedra, which we present in pseudocode followed by bounds on its
running time.  In particular, we show that the algorithm is
polynomial in the number of source images, when the dimension~$d$
is fixed.  (Later we shall state Conjecture~\ref{polynomial}, which
posits that the number of source images is polynomial in the number
of facets.)  Other applications, some of which are further
discussed in Section~\ref{limits}, include the discrete geodesic
problem (Corollary~\ref{c:geod}) and geodesic Voronoi diagrams
(Algorithm~\ref{voronoi-code} in Section~\ref{geod-voronoi}).

Roughly, Algorithm~\ref{code} consists of a single loop that with
every iteration constructs one new {\em event}.  Each event is a
pair consisting of a facet and a point that we have called a {\em
source image}\/ in the affine span of that facet.  The loop is
repeated exhaustively until all of the events are computed, so the
affine span of every facet has its full complement of source
images.  The Voronoi diagram for the set of source images in each
affine span induces a subdivision of the corresponding facet.  For
each maximal cell in this subdivision, the algorithm computes a
Euclidean motion (composition of rotation and parallel translation)
that moves it into the affine span of the facet
containing the source point.  The union of these moved images of
Voronoi cells is the output foldout~$\UUv$ in the tangent
space~$\Tv$ to the source point~$v$.

At each iteration of the loop, the algorithm must choose from a
number of {\em potential events}\/ that it could process into an
actual event.  Each potential event~$E$ consists of an already-computed
event $(\nu,F)$ plus a ridge~$R$ in the facet~$F$.  Processing the
event~$E$ applies a rotation to move the source image~$\nu$ into the
affine span of the other facet containing~$R$.  The potential event
that gets chosen must lie as close to the source point as possible;
this distance is the radius $r = r(R,\nu)$ at the beginning of the
loop.  The loop then calls Routine~\ref{r:event} to choose which
event to process; although this routine is quite simple in
structure, it is the part of the algorithm that most directly
encounters the subtlety of working in higher dimensions.  The end
of the loop consists of updating the sets of source images and
potential events; the latter requires Routine~\ref{r:cansee}, which
we have isolated because it is the only time-consuming part of the
algorithm, due to its Voronoi computation.

Let us emphasize that once a source point~$\nu$ is computed, it is
never removed.  This claim is part of Theorem~\ref{src}, in which
the correctness of Algorithm~\ref{code}---and indeed the procedure
of the algorithm itself---is more or less already implicit, as we
shall see in the proof of Theorem~\ref{correct}.

We assume that the convex polyhedron~$P$ is presented in the input
of the algorithm as an intersection of closed half-spaces.  Within
the algorithm, we omit the descriptions of standard geometric and
linear algebraic operations, for which we refer to~\cite{GO,Prep}.
These operations include the determination of lower-dimensional
faces (such as ridges) given the facets of~$P$, and the computation
of Voronoi diagrams.

Some additional notation will simplify our presentation of the
algorithm.  Denote by~$\cf$ and~$\CR$ the sets of facets and ridges
of~$P$, respectively.  If a ridge~\mbox{$R\in \CR$} lies in a facet
\mbox{$F\in \cf$}, denote by $\phi(F,R)$ the other facet
containing~$R$, so $F \cap \phi(F,R) = R$.  Finally, for each facet
$F \in \cf$, denote by $\EEE_F$ the set of all triples $(\nu,F,R)$
such that $\nu$ lies in the affine span~$\TF$ of~$F$, and $R \in
\CR$ is a ridge contained in~$F$.

\begin{alg} \label{code}
(Computing source unfolding)
\end{alg}
\begin{alglist}
\routine{input} convex polyhedron $P \subset \rr^{d+1}$ of
    dimension $d+1$\\
    point $v$ lying in the relative interior of a facet~$F$
    of~$P$
\routine{output} source foldout of the boundary $S = \partial P$
    into $\Tv \cong \rr^d$\, (see Sec.~\ref{unfold})
\begin{routinelist}{define}
\item[for each $F \in \cf$:] a finite set $\UF \subset \TF$ of points
\item[for each pair $(\nu,F)$] satisfying $\nu \in \UF$: an ordered
    list $\cl_{\nu,F}$ of facets
\item[for each $F \in \cf$:] a set $\EE_F \subset \EEE_F$ of
    \textsl{potential events}
\item[$\EE = \bigcup_{F \in \cf} \EE_F$, the set of \textsl{all
    potential events}]
\end{routinelist}
\begin{routinelist}{initialize}
\item[for $F \in \cf$: if $v \in F$ then $\UF := \nothing$ and
    $\EE_F = \nothing$;]
\item[otherwise $\UF := \{v\}$,] $\cl_{v,F} := (F)$,\ $\EE_F := \{(v,F,R)
    \mid R \in \CR$ and $R \subset F\}$
\end{routinelist}
\routine{compute} $\Phi_{F,F'}$ for all $F, F' \in \cf$ such that
    $F \cap F' \in \CR$ is a ridge\, (see Def.~\ref{fold})
\routine{while} $\EE \neq \nothing$
\begin{routinelist}{do}
    \item[$r := \min\{r(R,\nu) \mid (\nu,F,R) \in \EE\}$]\, (see
    Def.~\ref{cansee})
    \item[\procedure{choose a potential event $E = (\nu,F,R)
    \in \EE$ to process}]
    \item[set $F' := \phi(F,R)$, \
    $\nu' := \Phi_{F,F'}(\nu)$, \ $\cl_{\nu',F'} := (\cl_{\nu,F} \,F')$]
    \item[update $\UFp \gets \UFp \cup \{\nu'\}$]
    \item[\hskip1.35cm $\EE_{F'} \gets \{(\om,F',R') \in
        \EEE_{F'}$ such that $\om \in \UFp$, and]
    \item[\hskip2.82cm \procedure{point $\om \in \UFp$ can
        see~$R'$ through~$F'$}, and]
    \item[\hskip2.82cm $\om' \notin \UG$, where $G =
        \phi(F',R')$, $\om'= \Phi_{F',G}(\om)\}$]
    \item[\hskip1.35cm $\EE_F \gets \EE_F \minus \{E\}$,\
    $\EE\gets \bigcup_{G\in\cf} \EE_G$]
\end{routinelist}
\routine{end}{}\procedure{while-do}
\begin{routinelist}{compute}
    \item[for all facets $F \in \cf$ and points $\nu \in \UF$:]
    \item[\qquad\quad\ $\Phi_{\cl}$ for $\cl = \cl_{\nu,F}$\,
    (see Def.~\ref{sequential}), and then]
    \item[\qquad\quad\ $\UUv(\nu,F) := \Phi^{-1}_\cl\bigl(F\cap
    V(\UF,\nu)\bigr) \subset \Tv$]\, (see Thm.~\ref{vor})
\end{routinelist}
\routine{return} the foldout $\displaystyle\UUv = \bigcup_{(\nu,F)}
    \UUv(\nu,F)$, the union being over all $F \in \cf$ and $\nu
    \in \UF$
\end{alglist}
\vspace{-1ex}

\begin{rtne} \label{r:event}
(\procedure{choose a potential event to process})
\end{rtne}
\begin{alglist}
\routine{input} the set $\EE = \bigcup_{F \in \cf} \EE_F$ of
    potential events, and the radius~$r > 0$
\routine{output} an event $E \in \EE$\, (see Def.~\ref{ideal})
\routine{compute} the closest potential events
    $\EE_\circ :=\{(\om,F,R) \in \EE \mid r(\om,F) = r\}$
\routine{ } angle sequence $\angle(R,\om)$ for all
    $(\om,F,R) \in \EE_\circ$\, (see Def.~\ref{cansee})
\routine{find} a potential event $E = (\om,F,R) \in \EE_\circ$ with
    lexicographically
\routine{} minimal angle sequence $\angle(R,\om)$\,
    (see Sec.~\ref{wavefront})
\routine{return} the event $E = (\om,F,R)$
\end{alglist}

\begin{rtne} \label{r:cansee}
(\procedure{point $\om \in \Ups$ can see $R$ through $F$})
\end{rtne}
\begin{alglist}
\routine{input} facet $F \in \cf$, ridge $R \in \CR$, finite set of
    points $\Ups \ssu \TF$, and~$\om \in \Ups$
\routine{output} boolean variable $\beth \in\{\mathrm{True,False}\}$\,
    (see Def.~\ref{cansee})
\routine{compute}Voronoi diagram $\VV(\Ups)$\, (see Sec.~\ref{cutloci})
\routine{if} Voronoi cell $V(\Ups,\om) \ssu \VV(\Ups)$
    contains a point interior to~$R$\\
    and $\om$ lies on the same side of~$R$ as $F$ does in $\TF$\\
    then $\beth := \mathrm{True}$;\\
    otherwise $\beth := \mathrm{False}$
\routine{return} the variable $\beth$
\end{alglist}
\medskip

In the pseudocode, we have used the two different symbols `$\gets$'
and `$:=$' to distinguish between those variables that are being
updated and those that are being completely redefined at each
iteration of the \procedure{while-do} loop.  We hope this clarifies
the structure of Algorithm~\ref{code}.

\begin{thm} \label{correct}
For every convex polyhedron~$P \ssu \rr^{d+1}$ with boundary $S =
\partial P$, and any source point $v$ in a facet of~$S$,
Algorithm~\ref{code} computes the source foldout $\UUv \subseteq
\Tv$.
\end{thm}
\begin{proof}
First, we claim by induction that after each iteration of the
\procedure{while-do} loop, the set~$\{(\nu,F) \mid F\in \cf$ and $\nu
\in \UF\}$ is an order ideal in the source poset~$\src(v,S)$ from
Definition~\ref{poset}.  The claim is clear at the beginning of the
algorithm.  By construction, Routine~\ref{r:event} picks a minimal
potential event~$E$ to process.  The loop then adds an event by
processing~$E$, with the aid of Routine~\ref{r:cansee}.
Theorem~\ref{src} implies that what results after processing~$E$ is
still an order ideal of events, proving our claim.  Since the
poset~$\src(v,S)$ is finite by Lemma~\ref{finite}, the algorithm halts
after a finite number of loop iterations.  Finally, by
Theorem~\ref{vor} the Voronoi cells in each facet coincide with the
polyhedral subdivision of each facet by the cut locus~$\KKv$, so
Theorem~\ref{exp} shows that the foldout in the output is the desired
(nonoverlapping) {\em source foldout}\/~$\UUv$.%
\end{proof}

For purposes of complexity, we assume throughout this paper that
the dimension~$d$ is fixed.  Thus, if the convex polyhedron~$P \ssu
\rr^{d+1}$ of dimension~$d$ has $n$ facets, so $P$ is presented as
an intersection of~$n$ closed half-spaces, we can compute all of
the vertices and ridges of~$P$ in polynomial time~\cite{GO,Z}.  For
simplicity, we assume these are precomputed and appended to the
input.

The timing of Algorithm~\ref{code} crucially depends on the number
of source images.~~Let
\begin{eqnarray*}
  \upper &:=&  \max_{F \in \cf} \, |\src_F\bigr|
\end{eqnarray*}
be the largest number of source images in a tangent plane~$\TF$ for
a facet~$F$.  (This number can change if the source point~$v$ is
moved.  For example, $\upper = 4$ if $v$ is in the center of a
face, while $\upper = 12$ if $v$ is off-center as in
Figures~\ref{f:cube-voronoi} and~\ref{f:cube-source}.)  Note that
computing Voronoi diagrams for~$N$ points in~$\rr^d$ can be done in
$N^{O(d)}$ time \cite[p.381]{GO}.  See~\cite{Aur91,Cha91,For95} for
details and further references on Voronoi diagrams,
and~\cite{GO,Prep} for other geometric and linear \mbox{algebraic
computations we use}.

\begin{thm} \label{timing}
When the dimension~$d$ is fixed, the cost of Algorithm~\ref{code}
is polynomial in the number~$n$ of facets and the maximal
number\/~$\upper$ of source images for a facet.
\end{thm}
\begin{proof}
From the analysis in the proof of Theorem~\ref{correct}, the number of
loop iterations is at most $|\src(v,S)| \le |\cf|\upper\le n
\upper$.  Within the main body of the algorithm, only standard
geometric and linear algebraic operations are used, and these are all
polynomial in~$n$.  Similarly, Routine~\ref{r:event} uses only linear
algebraic operations for every potential event~$E \in \EE$.  Note that
the cardinality of the set of potential events~$\EE$ during any
iteration of the loop
is bounded by $|\src(v,S)| \cdot |\cf|^2 \le \bigl(n\upper\bigr)
\cdot n^2 = n^3 \upper$.

Routine~\ref{r:cansee} constructs Voronoi diagrams $\VV(\Ups)$ for
finite sets~$\Ups \ssu \rr^d$.  This computation requires
$|\Ups|^{O(d)} \le \bigl(\upper\bigr)^{O(d)}$ time, which is polynomial for our
fixed dimension~$d$.  Therefore the total cost of the algorithm is
also polynomial in~$n$ and~$\upper$.%
\end{proof}

\begin{cor} \label{c:geod}
Let $v$ and $w$ be two points on the boundary $S$ of the convex
$(d+1)$-dimensional polyhedron $P \ssu \rr^{d+1}$, and suppose
that~$v$ lies interior to a facet.  Then the geodesic
distance~$\mu(v,w)$ on~$S$ can be computed in time polynomial
in~$n$ and~$\upper$.%
\end{cor}

The restriction that $v$ lie interior to a facet is unnecessary,
and in fact Algorithm~\ref{code} can be made to work for arbitrary
points~$v$; see Sections~\ref{warpsource} and~\ref{geod-voronoi}.

\begin{proof} Use Algorithm~\ref{code} to compute the foldout map
$\vp: \UUv \to S$.  Find~$w' \in \Tv$ mapping to $w = \vp(w') \in S$,
and compute the distance~\mbox{$|v-w'|$}.  By the isometry of the
exponential map in Theorem~\ref{exp}, we conclude that $\mu(v,w) =
|v-w'|$.%
\end{proof}

\begin{remark} \label{r:lower}
The complexity of Algorithm~\ref{code} is exponential in~$d$ if the
dimension is allowed to grow.  For example, the number of vertices
of~$P$ can be as large as~$n^{\Omega(d)}$~\cite{Z}.
Similarly, the number of cells in
Voronoi diagrams of~$N$ points in~$\rr^d$ can be as large
as~$N^{\Omega(d)}$ \cite{Aur91,For95}.

On the other hand, for fixed dimension~$d$ Algorithm~\ref{code} can
not be substantially improved, because the input and the output
have costs bounded from below by (a~polynomial~in) $n$
and~$\upper$, respectively.  This is immediate for the input
since~$P$ is defined by~$n$ hyperplanes.  For the output, we claim
that the foldout~$\UUv$ in the output of Algorithm~\ref{code} can
not be presented at a smaller cost because it is a (usually
nonconvex) polyhedron that has at least $\upper$ {\em boundary
ridges}, meaning faces of dimension~\mbox{$d-1$} in the boundary
of~$\UUv$.  To see why, let $F$ be a facet with $\upper$ source
images, and for each $\nu \in \src_F$ consider a shortest path
$\ga_\nu$ whose sequential unfolding into $\TF$ has endpoint~$\nu$.
If instead we sequentially unfold the paths $\ga_\nu$ into~$\Tv$,
we get $|\src_F| = \upper$ segments emanating from~$v$.  Extend
each of these segments to an infinite ray.  Some of these infinite
rays might pierce the boundary of~$\UUv$ through faces of dimension
less than~\mbox{$d-1$}, but adjusting their directions slightly
ensures that each ray pierces the boundary of~$\UUv$ through a
boundary ridge.  These ridges are all distinct because their
corresponding rays traverse different facet sequences.

Of course, the efficiency of Algorithm~\ref{code} does not
necessarily imply that it yields an optimal solution to the
discrete geodesic problem---or the unfolding problem, for that
matter.  (The problem of computing {\em any}\/ nonoverlapping
unfolding, not necessarily the source unfolding, is of independent
interest in computational geometry~\cite{O}.)  But although
$\upper$ is not known to be polynomial in~$n$, we conjecture in
Section~\ref{open} that it is.  See Section~\ref{discretegeod} for
more history of the discrete geodesic problem.
\end{remark}

\begin{remark} \label{r:model}
Following traditions in computational geometry, we have not
specified our model of computation.  In most computational geometry
problems the model is actually irrelevant, since the algorithms are
oblivious to it.  In our case, however, the situation is more
delicate, due to the fact that during each iteration of the loop we
make a number of \emph{arithmetic operations}\/ that increase the
error.  More importantly, we make \emph{comparisons}, which
potentially require sharp precision.

Let us note that Theorem~\ref{timing} and its proof hold as stated
for the \emph{complexity over\/~$\rr$} model \cite{BCSS}, where
there are no errors, and all arithmetic operations and comparisons
have unit cost.  On the other hand, for the (usual)
\emph{complexity over\/~$\zz_2$} model \cite{BCSS}, the details are
less straightforward.  Think of the hyperplanes as given by
equations over the rational numbers~$\qqq$, and suppose that the
source point~$v$ is rational, as well.  Then the vertices, ridges,
source images, (tangents of angles in) angle sequences, and all
other data throughout the algorithm are also rational.  The whole
computation can be done \emph{precisely}\/ over~$\qqq$.  The
problem is keeping up with the computational precision as the
denominators grow exponentially.  Our preliminary calculations show
that Theorem~\ref{timing} still holds with several adjustments in
the proof.  Further exploration of this matter goes outside the
scope of this paper and~is~left~open.
\end{remark}

\unit{Convex polyhedral pseudomanifolds}\label{pseudomanifolds}

\noindent
Recall the notion of polyhedral complex from Section~\ref{geodesics}.
The results in Sections~\mbox{\ref{geodesics}--\ref{algorithm}} hold
with relatively little extra work for polyhedral complexes~$S$ that
are substantially more general than boundaries of polytopes.  Since
the generality is desirable from the point of view of topology, we
shall complete this extra work~here.

Suppose that~$x$ is a point in a polyhedral complex~$S$.  Denote by
\begin{eqnarray*}
  S_x(\ve) &=& \{y \in S \mid \mu(x,y) = \ve\}
\end{eqnarray*}
the {\em geodesic sphere}\/ in~$S$ at radius $\ve$ from~$x$.
If~$\<x\>$ is the smallest face of~$S$ containing~$x$, then for
sufficiently small positive real numbers~$\ve$, the intersection
$\<x\> \cap S_x(\ve)$ is an honest (Euclidean) sphere $\<x\>_\ve$ of
radius~$\ve$ around~$x$.  The set of points $N_x$ in~$S$ near~$x$ and
equidistant from all points on $\<x\>_\ve$ is the {\em normal space}\/
at~$x$ orthogonal to~$\<x\>$ in every face containing~$x$.  The {\em
spherical link}\/ of~$x$ at radius~$\ve$ is the set
\begin{eqnarray*}
  N_x(\ve) &=& \{y \in N_x \mid \mu(x,y) = \ve\}
\end{eqnarray*}
of points in the normal space at distance~$\ve$ from~$x$.  When $\ve$
is sufficiently small, the intersection of~$N_x(\ve)$ with any
$k$-dimensional face containing~$x$ is a sector inside a sphere of
dimension $k-1 - \dim\<x\>$.  The metric~$\mu$ on~$S$ induces a
subspace metric on the spherical link~$N_x(\ve)$.  Always assume $\ve$
is sufficiently small when $N_x(\ve)$ is written.

\begin{defn} \label{convex}
Let~$S$ be a connected finite polyhedral cell complex of dimension~$d$
whose facets all have dimension~$d$.  Given a point~$x$ inside the
union $S_{d-2}$ of all faces in~$S$ of dimension at most $d-2$, we say
that $S$~is {\em positively curved}\/ at~$x$ if the spherical
link~$N_x(\ve)$ is connected and has diameter less than~$\pi\ve$.  The
space $S$ is a {\em convex\/%
    \footnote{Using `convex' instead of `positively curved'
    allows usage of the term `nonconvex polyhedral complex'
    without ambiguity: `nonpositively curved' is already
    established in the context of CAT(0) spaces to mean (for
    polyhedral manifolds, at least) that no point has positive
    sectional curvature in any direction.  In contrast,
    `nonconvex' means that some point has a negative sectional
    curvature.}
polyhedral complex}\/ if $S$ is positively curved at every point $x
\in S_{d-2}$.
\end{defn}

This is definition of positive curvature is derived from the one
appearing in \cite{Sto76}.  It includes as special cases all
boundaries of convex polyhedra; this is essentially the content of
Proposition~\ref{warp}.



Spherical links give local information about geodesics, as noticed by
D.\thinspace{}A.~Stone (but see also \cite[Section~4.2.2]{BGP}).

\begin{lemma}[{\cite[Lemma~2.2]{Sto76}}] \label{join}
Suppose $S$ is a convex polyhedral complex.  Then $\tilde\ga$ is a
shortest path of length $\al\ve$ in the spherical link $N_x(\ve)$ of a
point $x \in S$ if and only if the union of all segments connecting
points of~$\tilde\ga$ to~$x$ is isometric (with distances given by the
metric on~$S$) to a sector of angle~$\al$ inside a disk in\/~$\rr^2$
of radius~$\ve$.
\end{lemma}

Although Stone only uses simplicial complexes, we omit the
straightforward generalization to polyhedral complexes.  Stone's
lemma forces shortest paths to avoid low-dimensional faces in the
presence of positive curvature.

\begin{prop} \label{convwarp}
Proposition~\ref{warp} holds for convex polyhedral complexes~$S$.
\end{prop}
\begin{proof}
Using notation from Lemma~\ref{join}, suppose that $\al < \pi$, and
let~$\ga$ be the segment connecting the endpoints of~$\tilde\ga$
through the sector of angle~$\al$.  Then $\ga$~misses~$x$.%
\end{proof}

The rest of Section~\ref{geodesics} goes through without change for
convex polyhedral complexes after we fix, once and for all, a {\em
tangent hyperplane}\/ $\TF \cong \rr^d$ for each facet~$F$.  The
choice of a tangent hyperplane is unique up to isometry.  For
convenience, we identify~$F$ with an isometric copy in~$\TF$, so that
(for instance) we may speak as if $F$ is contained inside~$\TF$.  This
makes Definition~\ref{fold}, in particular, work verbatim here.

The main difficulty to overcome in the remainder of
Sections~\ref{geodesics}--\ref{algorithm} is the finiteness in
Lemma~\ref{finite}.  In the context of convex polyhedral complexes,
this finiteness is fundamental.  It comes down to the fact that
shortest paths never wind arbitrarily many times around a single face
inside of a fixed small neighborhood of a point.  The statement of the
upcoming Proposition~\ref{fixed} would be false if we allowed
infinitely many facets, though it could still be made to hold in that
case if the sizes of the facets and their dihedral angles were forced
to be uniformly bounded away from zero.

\begin{prop} \label{fixed}
Fix a real number $r \geq 0$ and a convex polyhedral complex~$S$.
There is a fixed positive integer $N = N(r,S)$ such that the facet
sequence~$\cl_\ga$ of each shortest path~$\ga$ of length~$r$ in~$S$
has size~at~most~$N$.
\end{prop}
\begin{proof}
Pick a real number $\ve > 0$ small enough so that the following holds.
First, the sphere~$S_x(\ve)$ of radius~$\ve$ centered at each
vertex~$x$ only intersects faces containing~$x$.  Then, for every
point~$x$ on an edge but outside the union of the radius~$\ve$ balls
around vertices, the sphere~$S_x(\ve/2)$ only intersects faces
containing~$x$.  Iterating, for every point~$x$ on a face of
dimension~$i$ but outside the union of all the previously constructed
neighborhoods of smaller-dimensional faces, the sphere~$S_x(\ve/2^i)$
only intersects faces containing~$x$.  The existence of such a
number~$\ve$ follows from the fact that every facet of~$S$ is convex,
and that $S$ has finitely many facets (Definition~\ref{convex}).

It suffices to prove the lemma with \mbox{$r = \ve/2^d$}.  Let~$y$ be
the midpoint of~$\ga$.  The closed ball $B_y(\ve/2^{d+1})$ of
radius~$\ve/2^{d+1}$ centered at~$y$ intersects some collection of
faces, and among these there is a face of minimal dimension~$k$.  Fix
a point~$x_k$ lying in the intersection of this face with the ball
$B_y(\ve/2^{d+1})$.  The ball $B_{x_k}(\ve/2^k)$ contains~$\ga$ by the
triangle inequality.  However, $B_{x_k}(\ve/2^k)$ might also contain a
point~$x_j$ on a face of dimension $j < k$.  If so, then choose~$j$ to
be minimal.  Iterating this procedure (at most~$d$ times) eventually
results in a point~$x$ on face of dimension~$i$ such that
$B_x(\ve/2^i)$ contains~$\ga$ and only intersects faces
containing~$x$.

The metric geometry of~$S$ inside the ball $B_x(\ve/2^i)$ is the same
as in $B_{x'}(\ve/2^i)$ for every point~$x'$ on the smallest face
containing~$x$, as long as $B_{x'}(\ve/2^i)$ only intersects faces
containing~$x'$.  Since $S$ has finitely many faces by
Definition~\ref{convex}, we reduce to proving the lemma for shortest
paths~$\ga$ after replacing~$S$ by the ball $B = B_x(\ve/2^i)$.  In
fact, we shall uniformly bound the number of facets traversed by {\em
any}\/ shortest path~in~$B$.  For simplicity, inflate the metric by a
constant factor so that~$B$ has~radius~$2$.  By a {\em face of~$B$} we
mean the intersection of~$B$ with a face of~$S$.

Note that $B$ is isometric to a neighborhood of the apex on the
boundary of a right circular cone when the dimension is~\mbox{$d=2$}.
In this case, shortest paths in~$B$ can only pass at most once through
each ray emanating from~$x$.  We conclude that the lemma holds in full
(not just for~$B$) when $d=2$.  Using induction on~$d$, we shall
assume that the lemma holds in full for convex polyhedral complexes of
dimension at most~\mbox{$d-1$}.

First suppose that $x$ is not a vertex of~$S$, so the smallest
face~$\<x\>$ containing~$x$ has positive dimension.  Then $B$ is
isometric to a neighborhood of~$x$ in the product $\<x\> \times N_x$
of the face~$\<x\>$ with the normal space~$N_x$.  Projecting $\ga$
onto~$N_x$ yields a shortest path~$\bar\ga$ whose facet sequence in
the convex polyhedral complex~$N_x$ has the same size as~$\cl_\ga$.
Induction on~$d$ completes the proof in this case.
%
%

Now assume that~$x$ is a vertex of~$S$.  If one of the endpoints
of~$\ga$ is $x$ itself, then $\ga$ is contained in some face of~$B$.
Hence we may assume from now on that $x$ does not lie on~$\ga$.
Consider the radial projection from $B \minus \{x\}$ to the unit
sphere~$S_x(1)$ centered at~$x$ in~$B$.  If the image of~$\ga$ is a
point, then again $\ga$ lies in a single face; hence we may assume
that radial projection induces a bijection from~$\ga$ to its image
curve~$\tilde\ga$.  Since the geometry of~$B$ is scale invariant,
every path $\ga'$ in $B \minus \{x\}$ mapping bijectively
to~$\tilde\ga$ under radial projection has a well-defined facet
sequence equal to~$\cl_\ga$.

Choose another small real number~$\ve$ as in the first paragraph of
the proof, but with $B$ in place of~$S$.  Assume in addition that $\ve
< 1/2\pi$.  Subdivide~$\tilde\ga$ into at least $2^d/\ve$ equal arcs,
and use Lemma~\ref{join} to connect the endpoints of each arc by
straight segments in (the cone over~$\tilde\ga$ in)~$B$.
Lemma~\ref{join} implies that $\tilde\ga$ has length at most~$\pi$,
because $\ga$ is a shortest path.  Therefore each of the at least
$2^{d+1}\pi$ chords of~$\tilde\ga$ has length at most~$2^d$.  The
argument in the second paragraph of the proof now produces a new
center~$x'$ for each chord, and we are assured that $x' \neq x$
because the $\ve$-ball around~$x$ does not contain any of the chords.
Hence the smallest face~$\<x'\>$ containing~$x'$ has positive
dimension, and we are done by induction on~$d$ as before.%
\end{proof}

We shall see in Corollary~\ref{pseudofinite} that
Proposition~\ref{fixed} implies finiteness of the set of source
images.  But first, we need to introduce the class of polyhedral
complexes for which the notion of source image---and hence the rest of
Sections~\ref{geodesics}--\ref{algorithm}---makes sense.
%
%

\begin{defn}
A convex polyhedral complex~$S$ of dimension~$d$ is a {\em convex
polyhedral pseudomanifold}\/ if
$S$ satisfies two additional {\em pseudomanifold conditions\/}:
(i)~each facet is a bounded polytope of dimension~$d$, and (ii)~each
ridge lies in at most two facets.
\end{defn}

\begin{remark}
The ``A.D.\thinspace{}Aleks\-androv spaces with curvature bounded
below by~$0$'' of \cite{BGP} include convex polyhedral pseudomanifolds;
see Example~2.9(6) there.  Some of our results here, such as
surjectivity of exponential maps and non-branching of geodesics,
are general---and essentially local---properties of spaces with
curvature bounded below by~$0$.  But our focus is on decidedly
global issues pertaining to the combinatorial and polyhedral nature
of convex polyhedral pseudomanifolds, rather than on a local
analogy with Riemannian geometry.  That being said, many of our
results here can be extended to convex ``polyhedral''
pseudomanifolds with facets of constant positive curvature instead
of \mbox{curvature zero.  We leave this extension to the reader}.
\end{remark}

A flat point in an arbitrary convex polyhedral complex need not have a
neighborhood isometric to an open subset of~$\rr^d$, because more than
two facets could meet there.  In a convex polyhedral pseudomanifold,
on the other hand, every flat point not lying on the topological
boundary has a neighborhood isometric to an open subset of~$\rr^d$.
This condition is necessary for even the most basic of our results to
hold, including Corollary~\ref{cut} (whose proof works verbatim for
convex polyhedral pseudomanifolds), and the definition of source image
(which would require modification without it; see
Section~\ref{pseudo}).

We would have preferred to avoid the boundedness condition on facets,
but the finiteness of the set of source images in Lemma~\ref{finite}
can fail without it; see Section~\ref{finiteness}.

\begin{cor} \label{pseudofinite}
Lemma~\ref{finite} holds for convex polyhedral pseudomanifolds~$S$.
\end{cor}
\begin{proof}
Since every facet is bounded, the lengths of all shortest paths in~$S$
are uniformly bounded.  Proposition~\ref{fixed} therefore implies that
there are only finitely many possible facet sequences among all
shortest paths in~$S$ from the source.%
\end{proof}

Corollary~\ref{pseudofinite} yields the following consequences, with
the same proofs.

\begin{thm} \label{Svor}
Proposition~\ref{mount} on the generalization of Mount's lemma and
Theorem~\ref{vor} on Voronoi diagrams hold verbatim for convex
polyhedral pseudomanifolds~$S$.
\end{thm}

The rest of Section~\ref{cutloci} requires slight modification due to
the fact that a convex polyhedral pseudomanifold~$S$ can have a
nonempty topological boundary~$\partial S$.

\begin{prop} \label{pseudowarp-cut}
Fix a source point~$v$ in a convex polyhedral pseudomanifold~$S$.
Every warped point lies either in the topological boundary of~$S$ or
in the cut locus~$\KKv$.
\end{prop}
\begin{proof}
The same as Proposition~\ref{warp-cut}, assuming $w$ is not in the
boundary of~$S$.%
\end{proof}

In view of Proposition~\ref{pseudowarp-cut}, the statement of
Corollary~\ref{c:exp} fails for convex polyhedral pseudomanifolds.
Instead we get the following, with essentially the same proof.

\begin{cor} \label{c:pseudoexp}
If $v$ is a source point in a convex polyhedral pseudomanifold~$S$,
then
\begin{numbered}
\item \label{item:pseudopoly}
$\KKv \cup \partial S$ is a polyhedral complex of dimension $d-1$, and

\item \label{item:pseudoridge}
$\KKv \cup \partial S$ is the union $\Kv \cup S_{d-2} \cup \partial S$
of the cut, warped, and boundary points.
\end{numbered}
\end{cor}

The considerations in Section~\ref{unfold} go through with one small
modification: the noncompact flat Riemannian manifold $S^\circ$ is the
complement in~$S$ of not just the $(d-2)$-skeleton~$S_{d-2}$, but also
the topological boundary $\partial S$ of~$S$.  The notion of what it
means that a tangent vector at $w \in S$ {\em can be exponentiated}\/
(Definition~\ref{canbeexp}) remains unchanged, as long as $w$ lies
neither in $S_{d-2}$ nor the boundary of~$S$.  Similarly, the notion
of {\em source interior}\/ (Definition~\ref{U}) remains unchanged
except that the exponentials\/ \mbox{$\exp(t \zeta)$} for\/ $0 \leq t
\leq 1$ must lie in neither the cut locus~$\KKv$ nor the
boundary~$\partial S$.

\begin{thm} \label{pseudoexp}
Fix a source point~$v$ in the convex polyhedral pseudomanifold~$S$.
The exponential map\/ $\exp: \UUv \to S$ on the source foldout is a
polyhedral nonoverlapping foldout, and the boundary $\UUv \minus \Uv$
maps onto $\KKv \cup \partial S$.  Hence $\KKv \cup \partial S$ is a
cut set inducing a polyhedral nonoverlapping unfolding
\mbox{$S\!\minus\!(\KKv \cup \partial S) \to \Uv$ to the source
interior}.
\end{thm}
\begin{proof}
Using Corollary~\ref{c:pseudoexp} in place of Corollary~\ref{c:exp},
the proof is the same as that of Theorem~\ref{exp}, except that every
occurrence of $S \minus \KKv$ must be replaced by $S \minus (\KKv \cup
\partial S)$, and the open subspace $S^\circ$ must be defined as $S
\minus (S_{d-2} \cup \partial S)$ instead of $S \minus S_{d-2}$.%
\end{proof}

\begin{cor} \label{star-shaped}
Every convex polyhedral pseudomanifold of dimension~$d$ is, as a
metric space, obtained from a closed, star-shaped, polyhedral ball
in~$\rr^d$ by identifying pairs of isometric boundary components.
\end{cor}

Section~\ref{wavefront} concerns local geometry in the context of
convex polyhedra, and therefore requires no modification for convex
pseudomanifolds, given that all of the earlier results in the paper
hold in this more general context.

In Section~\ref{construct}, the only passage that does not seem to
work verbatim for convex polyhedral pseudomanifolds is the proof of
Proposition~\ref{jetsincrease}.  That proof is presented using
language as if~$F$ and~$F'$ were embedded in the same Euclidean
space~$\rr^{d+1}$, as they are in the case $S = \partial P$.  This
embedding can be arranged in the general case here by choosing
identifications of $\TF$ and $\TFp$ as subspaces of $\rr^{d+1}$ in
such a way that the copies of~$F$ and~$F'$ intersect as they do
in~$S$.

Finally, the algorithm in Section~\ref{algorithm} works just as well
for convex polyhedral pseudomanifolds, as long as these spaces are
presented in a manner that includes the structure of each facet as a
polytope and the adjacency relations among facets.  For example,
folding maps along ridges shared by adjacent facets can be represented
as linear transformations after assigning a vector space basis to each
tangent hyperplane.

For the record, let us summarize the previous three paragraphs.

\begin{thm}
The results in Sections~\ref{wavefront}, \ref{construct},
and~\ref{algorithm} hold verbatim for convex polyhedral
pseudomanifolds~$S$ in place of boundaries of convex polyhedra.
\end{thm}

\unit{Limitations, generalizations, and history}\label{limits}

\noindent
The main results in this paper are more or less sharp, in the sense
that further extension would make certain aspects of them false.  In
this section we make this sharpness precise, and also point out some
alternative generalizations of our results that might hold with
requisite modifications.  Along the way, we provide more history.

\subsection{Polyhedral vs.\ Riemannian} \label{wolter}

The study of geodesics on convex surfaces, where \mbox{$d=2$}, goes
back to ancient times and has been revived by Newton and the Bernoulli
brothers in modern times.  The study of explicit constructions of
geodesics on two dimensional polyhedral surfaces was initiated
in~\cite{L}, and is perhaps much older.

The idea of studying the exponential map on polyhedral surfaces goes
back to Aleksandrov~\cite[\S 9.5]{Al48}, who introduced it locally
when~$d=2$.  He referred to images of lines in the tangent space $\TF$
to a facet~$F$ as {\em quasi-geodesic}\/ lines on the surface, and
proved some results on them specific to the dimension~$d=2$.  Among
his other results was the $d=2$ case of Proposition~\ref{warp}.

A detailed analysis of the cut locus of $2$-dimensional convex
polyhedral surfaces was presented in~\cite{VP71}.  This paper,
seemingly overlooked in the West, gives a complete description of
certain convex regions called `peels' in \cite{AO}, which can be used
to construct source unfoldings.  The approach in~\cite{VP71} is
inherently $2$-dimensional and nonalgorithmic.

The study of exponential maps on Riemannian manifolds is
classical~\cite{Kob89}.  F.\thinspace{}Wolter \cite{Wol85} proved
properties of cut loci in the Riemannian context that are quite
similar to our results describing the cut locus as the closure of the
set of cut points.  In fact, we could deduce part~$2$ of our
Corollary~\ref{c:exp} from \cite[Lemma~2]{Wol85}---in the manifold
case, at least---using Proposition~\ref{warp} (which has no analogue
in Riemannian geometry).  The method would be to ``smooth out'' the
warped locus to make a sequence of complete Riemannian manifolds
converging (as metric spaces) to the polyhedral complex~$S$, such that
the complement of an ever decreasing neighborhood of the warped locus
in~$S$ is isometric to the corresponding subset in the approximating
manifold.  Every shortest path to $v$ in~$S$ is eventually contained
in the bulk complement of the smoothed neighborhood.

\begin{figure}[ht]
\begin{center}
\psfrag{v}{$v$}
\psfrag{w}{$w$}
\epsfig{file=\LOCAL/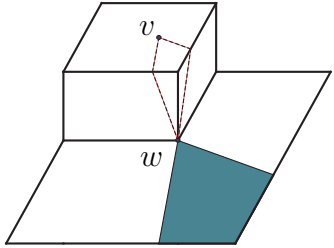,width=4.6cm}
\end{center}
\caption{Points in the shaded region have exactly two shortest
paths to~$v$; all of these paths go through the warped point~$w$.}
\label{f:ridgeset}
\end{figure}

This method does not extend to the polyhedral case where $S$ is
allowed to be nonconvex, because Proposition~\ref{warp} fails:
shortest paths (between flat points) can pass through warped points
(Fig.~\ref{f:ridgeset}).  Moreover, the polyhedrality in the first
part of Corollary~\ref{c:exp} fails systematically when~$S$ is
allowed to be nonconvex (see~\cite{MMP}).

\subsection{Low-dimensional flat faces}

We assumed in Definition~\ref{convex} that faces of
dimension~\mbox{$d-2$} or less in convex polyhedral complexes must be
nontrivially curved.  Allowing convex polyhedral complexes where
low-dimensional faces can be flat would break the notion of facet
sequence in Corollary~\ref{onepoint}, and would cause the set of
warped points to differ from the union of all closed faces of
dimension~\mbox{$d-2$}, in general.  The resulting definitions of
folding map and sequential unfolding would be cumbersome if not
completely opaque.  Nonetheless, the resulting definitions would be
possible, because shortest paths would still enter facets (and in
fact, all faces whose interiors are flat) at well-defined angles.  The
notion of exponential map would remain unchanged.

Definition~\ref{d:src} and Theorem~\ref{vor} should hold verbatim for
the modified notion of convex polyhedral pseudomanifold in which
low-dimensional flat faces are allowed, because the generalized Mount
Lemma (Proposition~\ref{mount}) should remain true.  Note that Mount's
lemma relies mainly on Proposition~\ref{warp} and Corollary~\ref{cut}.
The latter might be more difficult to verify in the presence of
low-dimensional flat faces, because it needs every flat point to have
a neighborhood isometric to an open subset of~$\rr^d$.  Thus one might
have to assume~$S$ is a {\em manifold}, and not just a pseudomanifold.

Observe that Fig.~\ref{f:mount} depends on not having low-dimensional
flat faces: it uses the fact that the vertex bordering the shaded
region must lie in the cut locus.

\subsection{Why the pseudomanifold conditions?} \label{pseudo}

Theorem~\ref{vor} fails for convex polyhedral complexes that are not
pseudomanifolds, even when there are no flat faces of small dimension.
Indeed, with the notion of cut point set forth in
Definition~\ref{locus},~entire facets could consist of cut points.  To
see why, suppose there is a cut point interior to a ridge lying on the
boundary of three or more facets, and note that the argument using
Fig.~\ref{f:Y} in the proof of Corollary~\ref{cut} fails.  For a more
concrete construction in dimension~\mbox{$d=2$}, find a convex
polyhedral pseudomanifold with a source point so that some edge in the
cut locus connects two vertices (for example, take a unit cube with a
source point in the center of a facet; see Fig.~\ref{f:two-fold}), and
then attach a triangle along that edge of the cut locus.  The attached
triangle (``dorsal fin'') \mbox{consists of cut points}.

The proof of Theorem~\ref{vor} fails for non-pseudomanifolds~$S$ when
we use the thinness of the cut set in the proof of
Proposition~\ref{mount}.  The appropriate definition of cut point~$x$
for convex polyhedral complexes more general than pseudomanifolds
should say that two shortest paths from~$x$ to the source leave~$x$ in
different directions---that is, they pierce the geodesic sphere
$S_x(\ve)$ at different points.  But Corollary~\ref{cut} would still
fail for shortest paths entering the ``dorsal fin'' constructed above.

\subsection{Aleksandrov unfoldings} \label{aleksandrov}

The dimension $d=2$ foldouts called `star unfoldings' in
\cite{AO,CH,AAOS} were conceived of by Aleksandrov in
\cite[p.~184]{Al48}.  Thus we propose here to use the term
`Aleksandrov unfolding' instead of `star unfolding', since in any case
these foldouts need not be star-shaped~polygons.  We remark that a
footnote on the same page in~\cite{Al48} indicates that Aleksandrov
did not realize the nonoverlapping property, which was only
established four decades later~\cite{AO}.

Aleksandrov unfoldings are defined for $3$-dimensional
polytopes~$P$ similarly to source unfoldings.  The idea is again to
fix a source point~$v$, but then slice the boundary~$S$ of~$P$ open
along each shortest path connecting~$v$ in~$S$ to a vertex.  An
example of the Aleksandrov unfolding of the cube is given in
\begin{figure}[ht]
\begin{center}
\epsfig{file=\LOCAL/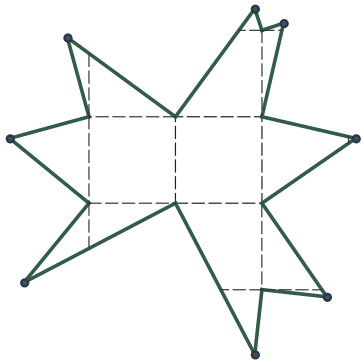,width=5cm}
\end{center}
\caption{Aleksandrov unfolding of the cube (the source point~$v$
is on the front face, while the left and back faces have no cuts).}
\label{f:cube-alex}
\end{figure}
Fig.~\ref{f:cube-alex} (see also Fig.~\ref{f:cube-top}).  Note that
when the source point is in the center of the face, the resulting
Aleksandrov unfolding agrees with the source unfolding in
Fig.~\ref{f:two-fold}.

\begin{figure}[ht]
\begin{center}
\epsfig{file=\LOCAL/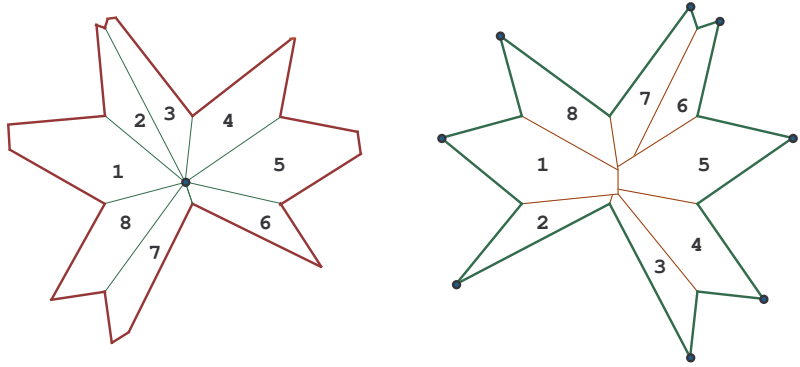,width=10.8cm}
\end{center}
\caption{Source and Aleksandrov unfoldings of the cube, where
the corresponding peels in both unfoldings are numbered from 1 to 8. }
\label{f:cube-connect}
\end{figure}

There is a formal connection between source and Aleksandrov
unfoldings.  Starting from the source unfolding, cut the star-shaped
polygon~$\UUv$ into sectors---these are `peels' as in
Section~\ref{wolter}---by slicing along the shortest paths to images
of vertices.  Rearranging the peels so that the various copies of~$v$
lie on the exterior cycle yields a nonoverlapping foldout \cite{AO}
containing an isometric copy of the bulk of the cut locus.  This
rearrangement is illustrated in Fig.~\ref{f:cube-connect}, which
continues Example~\ref{cube-source} (see~\cite{AAOS} for further
references).

\begin{figure}[ht]
\begin{center}
\epsfig{file=\LOCAL/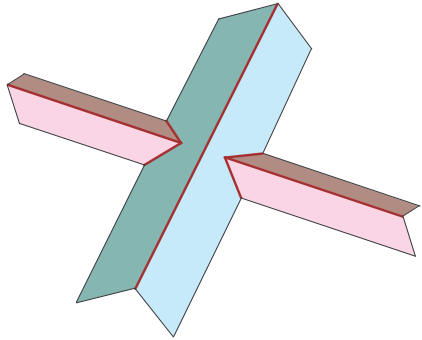,width=5.7cm}
\end{center}
\caption{Shortest paths to warped points making edges look
disconnected.}
\label{f:broken}
\end{figure}

No obvious higher-dimensional analogue of the Aleksandrov unfolding
exists, because although the union of all shortest paths connecting
the source point to warped points forms a polyhedral complex, this
complex is not a cut set as per Definition~\ref{d:unfold}.  Indeed,
thinking in terms of source foldouts again, the union of all rays
passing from the origin through the images of warped points does
not form the $(d-1)$-skeleton of a fan of polyhedral cones.  Even
when $d=3$, edges of~$S$ closer to the source point can make edges
farther away look disconnected, as seen from the source point.  An
example of how this phenomenon looks from~$v$ is illustrated in
Fig.~\ref{f:broken}, where the planes are meant to look tilted away
from the viewer.

Circumventing the above failure of Aleksandrov unfoldings in high
dimension would necessarily involve dealing with the fact that the
set~$S_{d-2}$ of warped points generically intersects the cut
locus~$\KKv$ in a polyhedral complex of dimension~\mbox{$d-3$}.
This ``warped cut locus'' usually contains points interior to
maximal faces of~$\KKv$, making it impossible for these interiors
of maximal cut faces to have neighborhoods in~$S$ isometric to open
sets in~$\rr^d$, even locally.  Thus the picture in
Fig.~\ref{f:cube-connect}, where most of the cut locus can lie
intact in~$\rr^2$, is impossible in dimension~\mbox{$d \geq 3$}.
The only remedy would be to make further slices across the
interiors of the maximal faces of the cut locus~$\KKv$ before
attempting to lay it flat in~$\rr^d$.  Making these extra slices in
a canonical way, to generalize Aleksandrov unfoldings to arbitrary
dimension, remains an open problem.

\subsection{Definition of source image}

Some subtle geometry dictated our choice of definition of `source
image' (Definition~\ref{d:src}).  With no extra information to go on,
we might alternatively have tried defining~$\src_F$ as the (finite)
set of endpoints of sequentially unfolded shortest~paths
\begin{itemize}
\item
ending at a point interior to~$F$; or

\item
ending anywhere on~$F$, including at a warped point.
\end{itemize}
Both look reasonable enough; but the first fails to detect faces of
dimension~\mbox{$d-1$} in the cut locus that lie entirely within
ridges of~$S$, while the second causes problems with verifying the
generalized Mount Lemma (Proposition~\ref{mount}) as well as
Proposition~\ref{jetsincrease} and Lemma~\ref{closest}.  Not that the
generalized Mount Lemma would be false with these ``bonus'' source
images included, but the already delicate proof would fail.
In~addition, having these extra source images would add unnecessary
bulk to the source poset.

\subsection{Finiteness of source images} \label{finiteness}

As we saw in Lemma~\ref{finite} for boundaries of polyhedra, or
Proposition~\ref{fixed} and Corollary~\ref{pseudofinite} for convex
polyhedral pseudomanifolds, the number of source images is finite.
The argument we gave in Lemma~\ref{finite} relies on the embedding
of $S$ as a polyhedral complex inside~$\rr^{d+1}$ in such a way
that each face is part of an affine subspace (i.e.\ not bent or
folded).  This embedding can be substituted by the more general
condition that the polyhedral metric on each facet is induced by
the metric on~$S$ (so pairs of points on a single facet are the
same distance apart in~$S$ as in the metric space consisting of the
isolated facet).  With this extra hypothesis, we would get
finiteness of the set of source images even for convex polyhedral
pseudomanifolds whose facets were allowed to be unbounded.
However, allowing unbounded facets in arbitrary convex polyhedral
pseudomanifolds can result in facets with infinitely many source
images.

\begin{figure}
\psfrag{S}{$S$}
\psfrag{U}{$\UUv$}
\psfrag{v}{\raisebox{-.5ex}{$v$}}
\epsfig{file=\LOCAL/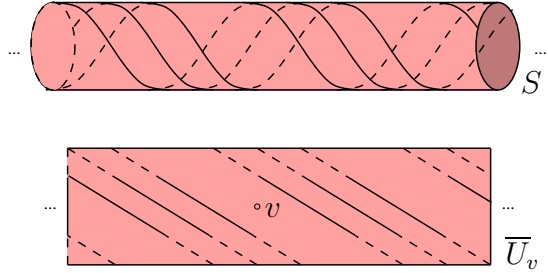,width=7.2cm}
\caption{A source foldout with infinitely many source images}
\label{f:cylinder}
\end{figure}

For example, consider an infinite strip in the plane, subdivided into
three substrips (one wide and two narrow, to make the picture
clearer).  Fix a distance $\ell > 0$, and glue each point on one
(infinite) boundary edge of the strip to the point $\ell$ units away
from its closest neighhbor on the opposite (infinite) edge of the
strip.  What results is the cylinder~$S$ in Fig.~\ref{f:cylinder}.
This cylinder would be a convex polyhedral manifold if its facets
were bounded.  The source foldout~$\UUv$ determined by a source
point~$v$ in the middle of the wide substrip is depicted
beneath~$S$, shrunken vertically by a factor of about~$2$.  The cut
locus~$\KKv$, which is a straight line along the spine of~$S$,
divides each substrip into infinitely many regions, so each
substrip has infinitely \mbox{many source images}.

\subsection{Generic source points} \label{generic}

For generic choices of source point~$v$, the source poset will be a
chain---that is, a total order on events.  The reason is that
moving~$v$ infinitesimally changes differently the angles from
different source images to ridges containing the same event point,
precisely because these source images are sequentially unfolded along
shortest paths leaving~$v$ in different directions.  Note, however,
that the {\em distances}\/ from these various source images to the
same event point always remain equal.


\subsection{Warped source points} \label{warpsource}

We assumed that the source point $v \in S$ lies in the relative
interior of some facet; however, nothing really changes when $v$ lies
in the relative interior of some ridge.  This can be seen by viewing
the exponential map as living on the interior flat points $S^\circ
\subseteq S$, as in Section~\ref{unfold}.

Moreover, simple modifications can generalize the exponential map to
the case where $v$ is warped.  However, exponentiation on the
complement of the cut locus cannot produce a nonoverlapping foldout
in~$\rr^d$ if $v \in S_{d-2}$, because the resulting cut locus would
not be a cut set.  Indeed, the cut locus would fail to contain all
of~$S_{d-2}$, so its complement could not possibly be isometric to an
open subset of~$\rr^d$.  On the other hand, exponentiation would
instead produce a foldout of~$S$ onto the tangent cone to~$S$ at~$v$.
The main point is that the source point still connects to a dense set
of points in~$S$ via shortest paths not passing through warped points,
by Proposition~\ref{warp}.

\subsection{Multiple source points} \label{geod-voronoi}

Let $\Ups = \{v_1.\ldots,v_k\} \ssu S$ be a finite set of points on
the boundary $S = \partial P$ of a convex polyhedron~$P$.  Define
the \emph{geodesic Voronoi diagram} $\VV_S(\Ups)$ to be the
subdivision of~$S$ whose closed cells are the sets
$$
V_S(\Ups,v_i) \, = \, \{w \in S  \mid \mu(v_i,w) \le \mu(v_j,w) \
\ \text{for all} \ \, 1 \le j \le k\}.
$$
Just like the (usual) Voronoi diagrams, computing geodesic Voronoi
diagrams is an important problem in computational geometry, with both
theoretical and practical applications \cite{AGSS,KWR,PL}
(see~\cite{Mit,GO} for additional references).

Below we modify Algorithm~\ref{code} to compute the geodesic Voronoi
diagrams in~$S$ when multiple source points are input.  The modified
algorithm outputs subdivisions of the facets of~$S$ that indicate
which source point is closest.  More importantly, it also computes
which combinatorial type of geodesic gives a shortest path.  In the
code below, the \procedure{while-do} loop and the routines remain
completely unchanged.  The only differences are in the initial and
final stages of the pseudocode.

\begin{alg} \label{voronoi-code}
(Computing the geodesic Voronoi diagram)
\end{alg}
\begin{alglist}
\routine{input} convex polyhedron $P \subset \rr^{d+1}$ of
    dimension $d+1$, and\\
    flat points $v_1,\ldots,v_k$ in the boundary $S = \partial P$
\routine{output} geodesic Voronoi diagram $\VV_S(\Ups)$ in~$S$

\routine{initialize} $\UF := \{v_i \mid v_i \in F\}$, \
    $\cl_{v_i,F} := (F)$ for $v_i \in \UF$, and\\
    \mbox{$\EE_F := \{(v_i,F,R) \in \EEE_F \mid$ \procedure{point
    $v_i \in \UF$ can see $R$ through $F\}$}}\\
    \textbf{[\,\ldots\,]}
\begin{routinelist}{compute}
    \item[for each~$i$: $\src_i := \{(\nu,F) \mid F \in \cf,
    \nu \in \UF$, and $\cl_{\nu,F}$ begins with $F_i\}$,] and
    \item[for each~$i$: the subset $V_S(\Ups,v_i) :=
    \bigcup_{(\nu,F)\in\src_i} V(\UF,\nu) \cap F$ of~$S$]
\end{routinelist}
\routine{return} geodesic Voronoi diagram
$\VV_S(\Ups) = \bigl(V_S(\Ups,v_1),\ldots,V_S(\Ups,v_k)\bigr)$
\end{alglist}
\bigskip

That some of the source points $v_1,\ldots,v_k$ might lie in the same
facet necessitates the call to Routine~\ref{r:cansee} in the
initialization of~$\EE_F$.  As we did before Theorem~\ref{timing},
define $\uppers$ to be the maximal number of
source images for a single facet.

\begin{thm} \label{voronoi-correct}
Let $P \ssu \rr^{d+1}$ be a convex polyhedron and $S = \partial P$,
with source points $v_1,\ldots,v_k$ in $S \minus S_{d-2}$.  For fixed
dimension~$d$, Algorithm~\ref{voronoi-code} computes the geodesic
Voronoi diagram~$\VV_S(\Ups)$ in time polynomial in~$k$, the
number~$n$ of facets, and~$\uppers$.
\end{thm}

The proof is a straightforward extension of the proof of
Theorems~\ref{correct} and~\ref{timing}; it is omitted.  Using
observations in Section~\ref{warpsource}, it is possible to modify
Algorithm~\ref{voronoi-code} to work for set of arbitrary (that is,
possibly warped) source points.

\subsection{The discrete geodesic problem} \label{discretegeod}

One of our motivating applications for this paper was to the {\em
discrete geodesic problem}\/ of computing geodesic distances and
the shortest paths between points $v$ and~$w$ in~$S$.  The
reduction of this problem to computing source unfoldings is easy:
construct the source foldout~$\UUv$ in the tangent cone at~$v$, and
compute the Euclidian distance between the images.

We should mention here that for $d=2$ essentially two methods are
used in the literature to resolve the discrete geodesic problem:
the construction of nonoverlapping unfoldings as above
(see~\cite{AAOS,CH,SS}), and the so-called `continuous Dijkstra'
method, generalizing Dijkstra's classical algorithm \cite{Dij59}
for finding shortest paths in graphs.  The second method originated
in~\cite{MMP} and is applicable to nonconvex surfaces (see
also~\cite{Ka}).  Interestingly, this method constructs an explicit
geodesic wavefront, and then selects and performs `events' one at a
time.  However, the time-ordering of events is based on the $d=2$
fact that the wavefront intersects the union of ridges (edges, in
this case) in a finite set of points.  Our approach is a
combination of these two algorithmic methods, which have previously
been separated in the literature.  We refer the reader to
\cite{Mit} for more references and results on the complexity of
discrete geodesic problems.  In general, computing geodesic
distances on arbitrary polyhedral complexes remains a challenging
problem of both theoretical and practical interest.

\unit{Open problems and complexity issues}\label{open}

\noindent
The source poset succeeds at time-ordering the events during wavefront
expansion, but it fails to describe accurately how the wavefront
bifurcates during expansion, because every event of radius less
than~$r$ occurs before the first event of radius~$r$ in the source
poset.  On the other hand, the notion of `geodesic precedence' from
Definition~\ref{geoprec} implies a combinatorial structure recording
bifurcation exactly.

\begin{defn} \label{tree}
Given a source point~$v$ on a convex polyhedral pseudomanifold~$S$,
the {\em vistal tree}\/ $\TT(v,S)$ is the set of events, partially
ordered by geodesic precedence.
\end{defn}

The definition of geodesic precedence immediately implies that
$\TT(v,S)$ is indeed a rooted tree.  It records the facet adjacency
graph of the polyhedral decomposition of the source foldout~$\UUv$
into cut cells of dimension~$d$.
Equivalently, this data describes the ``vista'' seen by an observer
located at the source point---that is, how the visual field of the
observer is locally subdivided by pieces of warped faces.
Proposition~\ref{jetsincrease} says precisely that the identity map on
the set of events induces a poset map from the vistal tree to the
source poset.  In particular, when the source point is generic as in
Section~\ref{generic}, the source poset is a linear extension of the
vistal tree.

There are numerous interesting questions to ask about the vistal tree,
owing to its geometric bearing on the nature of wavefront expansion on
convex polyhedra.  For example, its size, which is controlled by the
extent of branching at each node, is important for reasons of
computational complexity (Theorem~\ref{timing}).

\begin{conj} \label{polynomial}
The cardinality $|\src(v,S)|$ of the set of source images for a
polyhedral boundary~$S$ is polynomial in the number of facets when the
dimension~$d$ is~fixed.
\end{conj}

Hence we conjecture that there is a fixed polynomial~$f_d$,
independent of both~$S$ and~$v$, such that $|\src(v,S)| < f_d(n)$
for all boundaries $S = \partial P$ of convex polyhedra~$P$ of
dimension~\mbox{$d+1$} with $n$ facets, and all source points~$v
\in S$.  Note that the cardinality in question is at most factorial
in the number of facets: $|\src(v,S)| < n\cdot (n-1)! = n!$.
Indeed, each source image yields a facet sequence, and each of
these has length at most~$n$, starts at with facet~$F$
containing~$v$, and does not repeat any facet.
%
%

In fact, we believe a stronger statement than
Conjecture~\ref{polynomial} holds for boundaries of convex
polyhedra.  Given a shortest path~$\ga$, both of whose endpoints
lie interior to facets, call the facet sequence~$\LL_\ga$ traversed
by~$\ga$ the {\em combinatorial type}\/~of~$\ga$.

\begin{conj} \label{polynomial'}
The cardinality of the set of combinatorial types of shortest paths in
the boundary~$S$ of a convex polyhedron is polynomial in the number of
facets of~$S$, when the dimension is fixed.
\end{conj}

That is, we do not require one endpoint to be fixed at the source
point.  The result is stronger than Conjecture~\ref{polynomial}
because source images are in bijection with combinatorial types of
shortest paths in~$S$ with endpoint~$v$.  In both of the previous two
conjectures, the degree of the polynomial will increase with~$d$,
probably linearly.  When~$d=2$ both conjectures have been proved
(see~\cite{AAOS,CH}).

The intuition for Conjecture~\ref{polynomial} is that, as seen from
the source point in a convex polyhedral boundary~$S$, the faces of
dimension~\mbox{$d-2$} more or less subdivide the horizon into
regions.  (The horizon is simply the boundary of the source
foldout~$\UUv$, as seen from~$v$.)  The phrase `more or less' must be
made precise, of course; and our inability to delete it altogether is
a result of exactly the same phenomenon in Fig.~\ref{f:broken} that
breaks the notion of Aleksandrov unfoldings in higher dimension.

The reason we believe Conjecture~\ref{polynomial'} is that we believe
Conjecture~\ref{polynomial}, and there should not be too many
combinatorial types of vistal trees.  More precisely, moving the
source point a little bit should not alter the combinatorics of the
vistal tree, and there should not be more than polynomially many
possible vistal trees.  In fact, we believe a stronger, more geometric
statement.  It requires a new notion.

\begin{defn} \label{equivistal}
Two source points are {\em equivistal}\/ if their vistal trees are
isomorphic, and corresponding nodes represent the same facet
sequences.
\end{defn}

Again, the facet sequence corresponding to a node of the vistal tree
is the list of facets traversed by any shortest path whose sequential
unfolding yields the corresponding source image.  Hence two source
points are equivistal when their views of the horizon look
combinatorially the same.

\begin{conj} \label{vistal}
The equivalence relation induced by equivistality constitutes a convex
polyhedral subdivision of the boundary~$S$ of any convex polyhedron.
Moreover, the number of open regions in this subdivision is polynomial
in the number of\/ facets of\/~$S$.
\end{conj}

Independent from the conjecture's validity, the {\em vistal
subdivision}\/ it speaks of---whether convex polyhedral or not---is
{\em completely canonical}\/: it relies only on the metric
structure of~$S$.  In addition, lower-dimensional strata of the
vistal subdivision should reflect combinatorial transitions between
neighboring isomorphism classes of vistal trees.  Thus
Conjecture~\ref{vistal} gets at the heart of a number of issues
surrounding the interaction of the metric and combinatorial
structures of convex polyhedra.

\begin{remark}\label{r:geod-diam}
An important motivation behind the above ideas lies in the
computation of the geodesic diameter of the boundary of a convex
polytope.  This is a classical problem in computational geometry,
not unlike computing diameters of finite graphs (for the $d=2$ case
see~\cite{AO,AAOS}).  One possibility, for example, would be to
compute the vistal subdivision in Conjecture~\ref{vistal}, and use
this data to list the combinatorial types of shortest paths.  Each
combinatorial type could then be checked to determine how long its
corresponding shortest paths can be.  Conjectures~\ref{polynomial'}
and~\ref{vistal} give hope that the geodesic diameter problem can
be solved in polynomial~time.
\end{remark}

Let us remark here that the polynomial complexity conjectures fail
for nonconvex polyhedral manifolds of dimension~$d \geq 2$.  Note
that this does not contradict the fact that when~$d=2$ there exists
a polynomial time algorithm to solve the discrete geodesic problem
(see section~\ref{discretegeod} above).  Indeed, the number of
source images gives only a lower bound for \emph{our} algorithm,
while the problem is resolved by a different kind of algorithm.  On
the other hand, we show below that for $d \ge 3$ the discrete
geodesic problem is \textsc{NP}-hard.  The following result further
underscores the difference between the convex and nonconvex case.

\begin{prop} \label{nonconvex}
On (nonconvex) polyhedral manifolds, the number of combinatorial types
of shortest paths can be exponential in the number of facets.  In
addition, the problem of finding a shortest path on a (nonconvex)
polyhedral manifold is \textsc{NP}-hard.
\end{prop}

We present two proofs of the first part: one that is more explicit
and works for all~$d \ge 2$, and the other that is easy to modify
to prove the second part.  For the proof of the second part we
construct a $3$-dimensional polyhedral manifold, which is
essentially due to Canny and Reif~\cite{CR}.  See
Remark~\ref{r:nonconvex-surf} for comments on how to doctor these
manifolds to make them compact and without boundary.

\begin{proof}
To obtain a polyhedral domain with exponentially many shortest
paths between two points $x$ and~$y$, we consider a dimension $d=2$
example.  Simply take a pyramid shape polyhedral surface as shown
in Fig.~\ref{f:exp} and observe that there exist $2^k$ shortest
paths between top point~$v$ and bottom vertex~$w$, where~$k$ is the
number of terraces in the pyramid.  The omitted details are
straightforward.
\begin{figure}[hbt]
$$
\psfrag{v}{$v$}
\psfrag{w}{$w$}
\epsfig{file=\LOCAL/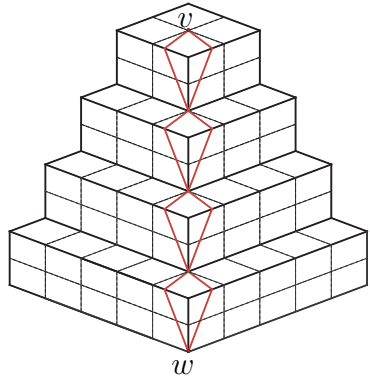,width=5cm}
$$
\caption{Nonconvex polyhedral surface in~$\rr^3$ and shortest
paths between points~$v$ and~$w$.}
\label{f:exp}
\end{figure}

Now consider a dimension $d=3$ example of a different type.
Polyhedrally subdivide~$\rr^3$ by taking the product of a
line~\mbox{$\ell = \rr$} with the subdivision of~$\rr^2$ in
Fig.~\ref{f:paths}.
\begin{figure}[hbt]
\epsfig{file=\LOCAL/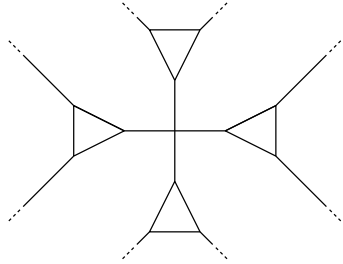,width=4.5cm}
\caption{Polyhedral subdivision of a planar slice in~\(\rr^3\)}
\label{f:paths}
\end{figure}
Observe that there are only finitely many cells.  Now add $4n$
hyperplanes $H_0,\ldots,H_{4n-1}$ orthogonal to~$\ell$, and equally
spaced along~$\ell$.  This still leaves finitely many convex cells.
Between hyperplanes $H_{4k}$ and~$H_{4k+1}$, for all $k=0 \dots n-1$,
remove all cells {\em
except}\/ the prisms whose bases are the top and bottom triangles in
Fig.~\ref{f:paths}.  Similarly, between hyperplanes
$H_{4k+2}$ and~$H_{4k+3}$, remove all cells {\em except}\/ the prisms
whose bases are the left
and right triangles in Fig.~\ref{f:paths}.

Now choose $x$ and~$y$ to be points on~$\ell$, with~$x$ being on one
side of all the hyperplanes, and~$y$ being on the other side.  Any
shortest path connecting $x$ to~$y$ must pass alternately through
vertical and horizontal pairs of triangular prisms, and there is no
preference for which of the two prisms in each pair the shortest path
chooses.  Thus the number of shortest paths is at least $4^n$, while
the number of cells is linear in~$n$.

For the second part, a construction in~\cite{CR} presents a
polyhedral domain~$B$ where the shortest path solution is
\textsc{NP}-hard.  This domain~$B$ is obtained by removing a set of
parallel equilateral triangles from~$\rr^3$.  To produce a manifold
one has to thicken the squares into nearly flat triangular
prisms. We omit the details.
\end{proof}

\begin{remark} \label{r:nonconvex-surf}
The polyhedral manifold~$S$ in the above proof is noncompact and
has nonempty boundary; but with a little extra work, we could
accomplish the same effect using a compact polyhedral manifold
without boundary.  The idea is to draw a large cube~$C$ around~$S$
in~$\rr^3$, and place copies $C_{\rm top}$ and~$C_{\rm bot}$ of~$C$
as the top and bottom facets of a hollow hypercube inside~$\rr^4$.
The remaining $6$ facets of the hollow hypercube are to remain
solid.  The result is compact, but still has nonempty boundary in
$C_{\rm top}$ and~$C_{\rm bot}$.  This we fix by building tall
$3$-dimensional prisms in~$\rr^4$ on the boundary faces, orthogonal
to $C_{\rm top}$ and~$C_{\rm bot}$, pointing away from the
hypercube.  Then we can cap off the prisms with copies of the cells
originally excised from $C \subset \rr^3$ to get a nonconvex
polyhedral $3$-sphere in~$\rr^4$.

%

\end{remark}

\begin{remark}\label{r:cayley}
The reader should not be surprised by the fact that computing the
geodesic distance is \textsc{NP}-hard for nonconvex manifolds.  On
the contrary: in most situations the problem of computing the
shortest distance is intractable, and in general is not
in~\textsc{NP}.  For example, finding the shortest distance in a
Cayley graph between two elements in a permutation group (presented
by a list of generators in~$S_N$) is known to be~\textsc{NP}-hard
even for abelian groups~\cite{EG}.  Furthermore, for directed
Cayley graphs the problem is \textsc{PSPACE}-complete~\cite{J}.
\end{remark}

% ...
%

Our final conjecture concerns the process of unfolding boundaries
of convex polyhedra: if someone provides a polyhedral
nonoverlapping foldout made of hinged wood, is it always possible
to glue its corresponding edges together?  Because wood is rigid,
we need not only a nonoverlapping property on the foldout as it
lies flat on the ground, but also a nonintersecting property as we
continuously fold it up to be glued.

Viewing this process in reverse, can we continuously unfold the
polyhedral boundary so that all dihedral angles monotonically
increase, until the whole polyhedral boundary lies flat on a
hyperplane?  This idea was inspired by recent works~\cite{BC,CDR} and
was suggested by Connelly.%
    \footnote{Private communication.}
While the monotone increase of the dihedral angles may seem an
unnecessary condition justified only by the aesthetics of the
blooming, it is in fact crucial in the references above.

As we have phrased things above, we asked for continuous unfolding of
an arbitrary nonoverlapping foldout.  But in fact, we only want to ask
that there {\em exist}\/ a foldout that can be continuously glued
without self-intersection.  Let us be more precise.

\begin{defn}
Let $S$ be the boundary of a convex polyhedron of dimension~$d+1$
in~$\rr^{d+1}$.  A~{\em continuous blooming}\/ of~$S$ is a choice of
nonoverlapping foldout $\ol U \to S$, and a homotopy $\{\phi_t : \ol U
\to \rr^{d+1} \mid 0 \leq t \leq 1\}$ such that
\begin{numbered}
\item
$\phi_0$ is the foldout map $\ol U \to S$;

\item
$\phi_1$ is the identity map on~$\ol U$;

\item
$\phi_t$ is an isometry from the interior $U$ of~$\ol U$ to its image,
and $\phi_t$ is linear on each component of the complement of the cut
set in each facet, for $0 < t < 1$; and

\item
the dihedral angles between corresponding facets of $\phi_t(\ol U)$
increase as $t$~increases.
\end{numbered}
\end{defn}

An example of a continuous blooming is given in Figure~\ref{f:bloom} below.

\begin{conj} \label{bloom}
Every convex polyhedral boundary has a continuous blooming.
\end{conj}

\begin{figure}[hbt]
$$
\epsfig{file=\LOCAL/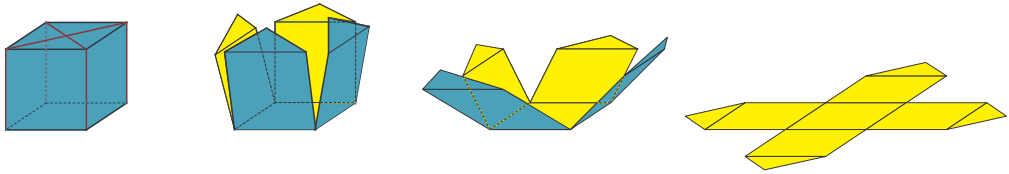,width=13.5cm}
$$
\caption{An example of a continuous blooming of the surface of the
cube.}
\label{f:bloom}
\end{figure}

Even though we ask only for existence, we believe that in fact the
source unfolding can be continuously bloomed.  As far as we know, this
is open even for $d=2$.  Interestingly, we know of no nonoverlapping
unfolding that cannot be continuously bloomed, and remain in
disagreement on their potential existence.

\subsection*{Acknowledgments} \label{ackpage}

We express our gratitude to Laci Babai, Yuri Burago, Bob Connelly,
Erik Demaine, Martin Demaine, Maksym Fedorchuk, Robin Forman, Sam
Grushevsky, Tracy Hall, Bob MacPherson, Jon McCammond, David Mount,
J\'anos Pach, Grigori Perel'man, Micha Sharir, Victor Shnayder, and
Santosh Vempala for helpful conversations.  We are indebted to
Robin Forman, Galina Shubina, and Frank Wolter for providing
references, and to Frank Sottile for pointing out a subtle error in
an earlier version.
The second author drew inspiration from~\cite{Dali}.

The first author was partially supported by the National Science
Foundation, and he wishes to acknowledge the Mathematisches
Forschungsinstitut Oberwolfach for a stimulating research
environment during the week-long program on Topological and
Geometric Combinatorics (April, 2003).  Most of this work was
completed while the first author was at the Massachusetts Institute
of Technology (Cambrdige, MA) and the Mathematical Sciences
Research Institute (Berkeley, CA).  The second author was partially
supported by the National Security Agency and the National Science
Foundation.  He would like to thank the organizers of the ``Second
Geometry Meeting dedicated to A.D.\thinspace{}Aleksandrov'' held at
the Euler International Mathematical Institute in St. Petersburg
(June, 2002), where these results were originally presented.


\begin{thebibliography}{BDEKMS}\label{refpage}

\bibitem[AAOS97]{AAOS}
P.~K. Agarwal, B. Aronov, J. O'Rourke, and
  C.~A. Schevon, Star unfolding of a polytope with applications,
  \emph{SIAM J. Comput.} \textbf{26} (1997), no.~6, 1689--1713.

\bibitem[AGSS89]{AGSS}
A. Aggarwal, L.~J. Guibas, J. Saxe, and  P.~W. Shor,
  A linear-time algorithm for computing the Voronoi diagram of a
  convex polygon, \emph{Discrete \& Comput. Geom.} \textbf{4}
  (1989), no.~6, 591--604.


\bibitem[Ale48]{Al48}
A.~D. Aleksandrov, \emph{Vnutrennyaya geometriya vypuklykh
  poverkhnostey {\rm (in Russian)}}, {\mbox{M.--L.}: Gostekhizdat},
  1948.

\bibitem[Ale50]{Al50}
A.~D. Aleksandrov, \emph{Vypuklye mnogogranniki {\rm (in
  Russian)}}, {M.: Gostekhizdat}, 1950.

\bibitem[AO92]{AO}
B. Aronov and J. O'Rourke, Nonoverlap of the star unfolding,
  \emph{Discrete \& Comput. Geom.} \textbf{8} (1992), no.~3,
  219--250.

\bibitem[Aur91]{Aur91}
F. Aurenhammer, Voronoi diagrams---a survey of a fundamental
  geometric data structure, \emph{ACM Comput. Surv.} \textbf{23}
  (1991), 345--405.

\bibitem[BC02]{BC}
K. Bezdek and R. Connelly,  Pushing disks apart---the
  Kneser--Poulsen conjecture in the plane, \emph{J. Reine
  Angew. Math.} \textbf{553}  (2002), 221--236.

\bibitem[BDEKMS]{BDEK}
M. Bern, E.~D. Demaine, D. Eppstein, E. Kuo, A. Mantler, and
  J. Snoeyink, Ununfoldable polyhedra with convex faces,
  \emph{Comput. Geom.} \textbf{24} (2003),  51--62.

\bibitem[BCSS98]{BCSS}
L. Blum, F. Cucker, M. Shub, and S. Smale, \emph{Complexity and
  real computation}, Springer, New York, 1998.

\bibitem[BGP92]{BGP}
Yu. Burago, M. Gromov, and G. Perelman, A.D.\thinspace{}Alexandrov
  spaces with curvature bounded below, \emph{Russian Math. Surveys}
  \textbf{47} (1992), no.~2, 1--58.

\bibitem[CR87]{CR}
J.~F. Canny and J.~H. Reif, New lower bound techniques for robot
  motion planning problems, \emph{Proc. 28-th IEEE FOCS} (1987),
  49--60.

\bibitem[Cha91]{Cha91}
B. Chazelle, An optimal convex hull algorithm and new results on
  cuttings, \emph{Proc. 32-nd IEEE FOCS} (1991), 29--38.

\bibitem[CH96]{CH}
J. Chen and Y. Han, Shortest paths on a polyhedron. I. Computing
  shortest paths, \emph{Internat. J. Comput. Geom. Appl.} {\bf 6}
  (1996), no. 2, 127--144.

\bibitem[CDR03]{CDR}
R. Connelly, E.~D. Demaine and G. Rote, Straightening polygonal
  arcs and convexifying polygonal cycles, \emph{Discrete \&
  Comp. Geom.}, to appear.

\bibitem[Dal54]{Dali}
S. Dali, {\em The Crucifixion} (`\emph{Corpus Hypercubus}'),
  Lithograph, Port Lligat, 1954, owned by  Humboldt State
  University, Arcata, CA.

\bibitem[Dij59]{Dij59}
E.~W. Dijkstra, A note on two problems in connexion with graphs,
  \emph{Numer. Math.} \textbf{1} (1959), 269--271.

\bibitem[EG81]{EG}
S. Even and  O. Goldreich, The minimum-length generator sequence
  problem is NP-hard, \emph{J. Algorithms} \textbf{2} (1981),
  no. 3, 311--313.

\bibitem[For95]{For95}
S. Fortune, Voronoi diagrams and {D}elaunay triangulations, in
  \emph{Computing in Euclidean geometry} (F.~Hwang and D.~Z. Du,
  eds.), World Scientific, Singapore, 1995, 225--265.

\bibitem[GO97]{GO}
J.~E. Goodman and J. O'Rourke (eds.), \emph{Handbook of discrete
  and computational geometry}, CRC Press, Boca Raton, FL, 1997.

\bibitem[Jer85]{J}
M.~R. Jerrum, The complexity of finding minimum-length generator
  sequences, \emph{Theoret. Comput. Sci.} \textbf{36} (1985),
  no. 2-3, 265--289

\bibitem[Kap99]{Ka}
S. Kapoor, An efficient computation of geodesic shortest paths,
  \emph{Proc. of the 31-st ACM STOC} (1999), 770--779.

\bibitem[Kob89]{Kob89}
S. Kobayashi, On conjugate and cut loci, in {\em Global
  differential geometry}, MAA, Washington, DC, 1989, pp.~140--169.

\bibitem[KWR97]{KWR}
R. Kunze, F.~E. Wolter, T. Rausch, Geodesic Voronoi Diagrams on
  Parametric Surfaces, in \emph{Proc. Comp. Graphics Int.} (1997),
  Hasselt-Diepenbeek, Belgium, 230--237.

\bibitem[Lyu40]{L}
L.~A. Lyusternik, \emph{Geodesic lines.  The shortest paths on
  surfaces {\rm (in Russian)}}, {\mbox{M.--L.}: Gostekhizdat}, 1940.

\bibitem[Mit00]{Mit}
J.~S.~B. Mitchell, Geometric shortest paths and network
  optimization, in \emph{Handbook of computational geometry},
  633--701, North-Holland, Amsterdam, 2000.

\bibitem[MMP87]{MMP}
J.~S.~B. Mitchell, D.~M. Mount, and C.~H. Papadimitriou,  The
  discrete geodesic problem, \emph{SIAM J.  Comp.} \textbf{16}
  (1987), no.~4, 647--668.


\bibitem[Mou85]{Mount85}
D.~M. Mount, \emph{On finding shortest paths on convex polyhedra},
  Technical Report 1495, Dept. of Computer Science, Univ. of
  Maryland, Baltimore, MD, 1985.

\bibitem[O'R00]{O}
J. O'Rourke, Folding and unfolding in computational geometry, in
  \emph{Discrete and Computational Geometry} (Tokyo, 1998),
  Springer, Berlin, 2000, pp.~258--266.


\bibitem[PL98]{PL}
E. Papadopoulou and D.~T. Lee,  A new approach for the geodesic
  Voronoi diagram of points in a simple polygon and other
  restricted polygonal domains, \emph{Algorithmica} \textbf{20}
  (1998), no.~4, 319--352.

\bibitem[PS85]{Prep}
F.~P. Preparata and M.~I. Shamos, \emph{Computational geometry. An
  introduction}, Texts and Monographs in Computer Science,
  Springer, New York, 1985.


\bibitem[SS86]{SS}
M. Sharir and A. Schorr, On shortest paths in polyhedral spaces,
  \emph{SIAM J. Comp.} \textbf{15} (1986), no.~1, 193--215.

\bibitem[Sto76]{Sto76}
D.~A. Stone, Geodesics in piecewise linear manifolds,
  \emph{Trans. Amer. Math. Soc.} \textbf{215} (1976), 1--44.

\bibitem[VP71]{VP71}
Ju.~A. Volkov and E.~G. Podgornova, The cut locus of a polyhedral
  surface of positive curvature (in Russian), \emph{Ukrainian
  Geometric Sbornik} \textbf{11} (1971), 15--25.

\bibitem[Wol85]{Wol85}
F.~E. Wolter, \emph{Cut loci in bordered and unbordered Riemannian
  manifolds}, Ph.D. thesis, TU Berlin, FB Mathematik, Berlin,
  Germany, 1985.

\bibitem[Zie95]{Z}
G.~M. Ziegler, \emph{Lectures on polytopes}, Springer, New York,
  1995.

\end{thebibliography}

\end{document}